\documentclass[smallextended,final,natbib]{svjour3}
\usepackage[]{graphicx}
\usepackage{subeqn,url}
\usepackage{times,bbold}
\newcommand{\DDD}{^{\mathrm{3D}}}
\newcommand{\DD}{^{\mathrm{2D}}}
\newcommand{\Dkkp}{D(\bk,\bk')}
\newcommand{\Dkpkc}{D^*(\bk',\bk)}
\newcommand{\Dlmlmp}{D_{lm,l'm'}} 

\newcommand{\Drhrhp}{D(\rhat,\rhat')} 
\newcommand{\Dttp}{D(t,t')}
\newcommand{\Dxixip}{D(\bxi,\bxi')}
\newcommand{\Dxixi}{D(\bxi,\bxi)}
\newcommand{\Dxxp}{D(\bx,\bx')}
\newcommand{\Dxx}{D(\bx,\bx)}
\newcommand{\D}{^{\mathrm{1D}}}
\newcommand{\Lpot}{(L+1)^2}

\newcommand{\Trm}{^{\mathrm{\sst{T}}}}

 \newcommand{\Ylmp}{Y_{l'm'}}
\newcommand{\Ylmrh}{Y_{lm}(\rhat)}
\newcommand{\Ylm}{Y_{lm}}
\newcommand{\also}{\qquad\mbox{and}\qquad} 
\newcommand{\barray}{\begin{array}}

\newcommand{\ber}{\begin{eqnarray}} 
\newcommand{\be}{\begin{equation}}

\newcommand{\bk}{\mathbf{k}}

\newcommand{\bxi}{\mbox{\boldmath$\xi$}}
\newcommand{\bx}{\mathbf{x}}
\newcommand{\sfg}{\mathsf{g}}
\newcommand{\sfh}{\mathsf{h}}
\newcommand{\dab}{\delta_{\alpha\beta}}
\newcommand{\dbk}{\,d\bk}
\newcommand{\dbx}{\,d\bx}
\newcommand{\dk}{\,d k}

\newcommand{\domg}{\,d\Omega}

\newcommand{\dth}{\,d \theta}
\newcommand{\earray}{\end{array}}
\newcommand{\edc}{\end{document}}
\newcommand{\eer}{\end{eqnarray}}
\newcommand{\ee}{\end{equation}}
\newcommand{\ena}{\end{eqnarray}}
\newcommand{\eqa}{\begin{eqnarray}}
\newcommand{\fnorm}{(2\pi)^{-2}}

\newcommand{\fracd}[2]{\frac{\displaystyle{#1}}{\displaystyle{#2}}}

\newcommand{\glmp}{g_{l'm'}} 
\newcommand{\glm}{g_{lm}} 
\newcommand{\grh}{g(\rhat)}
\newcommand{\grhp}{g(\rhat')} 
\newcommand{\half}{\frac{1}{2}}
\newcommand{\hspo}{\hspace{0.1em}}
\newcommand{\hspf}{\hspace{0.05em}}

\newcommand{\hsp}{\hspace*{0.1em}}
\newcommand{\intK}{\int_\mathcal{K}}
\newcommand{\intR}{\int_\mathcal{R}}
\newcommand{\intW}{\int_{-W}^{W}}

\newcommand{\intinf}{\int_{-\infty}^{\infty}}
\newcommand{\intinft}{\int_{\mathbb{R}^2}}
\newcommand{\intk}{\int_0^K}
\newcommand{\into}{\int_\Omega}

\newcommand{\intr}{\int_R} 
\newcommand{\intth}{\int_0^{2\pi}}
\newcommand{\intzinf}{\int_{0}^{\infty}}
\newcommand{\jnorm}{(2\pi)^{-1}}

\newcommand{\mbf}{\mathbf}

\newcommand{\nnr}{\nonumber}

\newcommand{\rar}{\rightarrow}
\newcommand{\rhat}{\mbf{\hat{r}}} 
\newcommand{\ssec}{\subsection}

\newcommand{\ssK}{\mathcal{K}}
\newcommand{\ssL}{\mathcal{L}} 
\newcommand{\ssP}{\mathcal{P}}
\newcommand{\ssQ}{\mathcal{Q}}
\newcommand{\ssQi}{\mathcal{Q}^{-1}\!}
\newcommand{\ssR}{\mathcal{R}}

\newcommand{\sst}{\scriptstyle} 
\newcommand{\suml}{\sum\limits}
\newcommand{\sumshLp}{\suml_{l'=0}^{L}\suml_{m'=-l'}^{l'}}
\newcommand{\sumshL}{\suml_{l=0}^{L}\suml_{m=-l}^{l}}
\newcommand{\sumsh}{\suml_{l=0}^{\infty}\suml_{m=-l}^{l}}
\newcommand{\tlofp}{\left(\frac{2l+1}{4\pi}\right)}
 
\newcommand{\with}{\quad\mbox{with}\quad} 
\begin{document}
\title{Spatiospectral concentration in the Cartesian plane}
\author{Frederik J.~Simons\and Dong
  V.~Wang\thanks{Now at: Department of Statistics and 
    Operations Research, The University of North Carolina at Chapel
    Hill, Chapel Hill, NC 27599, U.S.A.}}
\institute{Department of Geosciences, Princeton University, Princeton
  NJ 08544, U.S.A.}
\maketitle

\begin{abstract}
We pose and solve the analogue of Slepian's time-frequency
concentration problem in the two-dimensional plane, for applications
in the natural sciences. We determine an orthogonal family of strictly
bandlimited functions that are optimally concentrated within a closed
region of the plane, or, alternatively, of strictly spacelimited
functions that are optimally concentrated in the Fourier domain. The
Cartesian Slepian functions can be found by solving a Fredholm
integral equation whose associated eigenvalues are a measure of the
spatiospectral concentration. Both the spatial and spectral regions of
concentration can, in principle, have arbitrary geometry. However, for
practical applications of signal representation or spectral analysis
such as exist in geophysics or astronomy, in physical space irregular
shapes, and in spectral space symmetric domains will usually be
preferred. When the concentration domains are circularly symmetric in
both spaces, the Slepian functions are also eigenfunctions of a
Sturm--Liouville operator, leading to special algorithms for this case,
as is well known. Much like their one-dimensional and spherical
counterparts with which we discuss them in a common framework, a basis
of functions that are simultaneously spatially and spectrally
localized on arbitrary Cartesian domains will be of great utility in
many scientific disciplines, but especially in the geosciences. 
\end{abstract}


\section{Introduction}

The one-dimensional prolate spheroidal wave functions (pswf) have
enjoyed an enduring popularity in the signal processing community ever
since their introduction in the
early~1960s~\cite[]{Landau+61,Landau+62,Slepian+61}. Indeed, in many
scientific and engineering disciplines the pswf and their relatives
the discrete prolate spheroidal sequences
(dpss)~\cite[]{Grunbaum81a,Slepian78} have by now become the preferred
data windows to regularize the quadratic inverse problem of power
spectral estimation from time-series observations of finite
extent~\cite[]{Percival+93}. At the deliberate cost of introducing
spectral bias, windowing the data with a set of such orthogonal
``tapers'' lowers the variance of the ``multitaper''
average~\cite[]{Thomson82}, which results in estimates of the power
spectral density that are low-error in the mean-squared
sense~\cite[]{Thomson90}. As a basis for function representation,
approximation and
interpolation~\cite[]{Delsarte+85,Moore+2004,Shkolnisky+2006,Xiao+2001},
or in stochastic linear inverse
problems~\cite[]{Bertero+85,Bertero+88,DeVilliers+2001,Wingham92}, the
pswf have been less in the public eye, especially compared to wavelet
analysis~\cite[]{Daubechies92,Percival+2006}, though, due to advances in
computation, there has been a resurgent interest in recent
years~\cite[]{Beylkin+2002,Karoui+2008,Khare+2003,Walter+2005b}, in
particular as relates to using them for the numerical solution of
partial differential
equations~\cite[]{Beylkin+2005,Boyd2003,Boyd2004,Chen+2005}.

The pswf are the solutions to what has come to be known as ``the
concentration problem''~\cite[]{Flandrin98,Percival+93} of Slepian,
Landau and Pollak, in which the energy of a bandlimited function is
maximized by quadratic optimization inside a certain interval of
time. Vice versa, it refers to the maximization of the spectral
localization of a timelimited function inside a certain target
bandwidth~\cite[]{Slepian83}. In the first version of the problem, the
bandlimited, time-concentrated pswf form an orthogonal basis for the
entire space of bandlimited signals that is also orthogonal over the
particular time interval of interest. In the second version the
timelimited, band-concentrated pswf are a basis for square-integrable
broadband signals that are exactly confined to the
interval~\cite[]{Landau+61,Slepian+61}. In general we shall refer to all
singular functions of time-bandwidth or space-bandwidth projection
operators as ``Slepian functions''.

The fixed prescription of the ``region of interest'' in physical or
spectral space is a deliberately narrow point of view that is well
suited to scientific or engineering studies where the assumption of
stationarity, prior information, or the availability of data will
dictate the interval of study from the outset. This distinguishes the
Slepian functions philosophically from the eigenfunctions of
full-phase-space localization
operators~\cite[]{Daubechies88a,Simons+2003a} or
wavelets~\cite[]{Daubechies+88,Olhede+2002}, with which they
nevertheless share strong
connections~\cite[]{Lilly+95,Shepp+2000,Walter+2004,Walter+2005a}. Strict
localization of this type remained the driving force behind the
development of Slepian functions over fixed geographical domains on
the surface of the
sphere~\cite[]{Albertella+99,Miranian2004,Simons+2006a} --- which have
numerous applications in
geodesy~\cite[]{Albertella+2001,Han+2008b,Simons+2006b},
geomagnetism~\cite[]{Simons+2009b,Schott+2011},
geophysics~\cite[]{Han+2007,Han+2008c,Han+2008a,Harig+2010},
biomedical~\cite[]{Maniar+2005,Mitra+2006} and
planetary~\cite[]{Evans+2010,Han2008,Han+2009,Wieczorek+2005} science, and
cosmology~\cite[]{Dahlen+2008,Wieczorek+2007} --- as opposed to
approaches using spherical
wavelets~\cite[]{Chambodut+2005,Fay+2008,Fengler+2006a,Freeden+97,Holschneider+2003,Kido+2003,McEwen+2007,Panet+2006,Schmidt+2006},
needlets~\cite[]{Marinucci+2008},
splines~\cite[]{Amirbekyan+2008b,Lai+2009,Michel+2008}, radial basis
functions~\cite[]{Freeden+99,Schmidt+2007}, coherent
states~\cite[]{Hall+2002,Kowalski+2000,Tegmark95,Tegmark96a}, or other
constructions~\cite[]{Simons+97a}, which have all come of age in these
fields also. Finally, we note that the very choice of the criterion to
define localization is open to
discussion~\cite[]{Narcowich+96,Parks+90,Riedel+95,Saito2007,Wei+2010};
in particular, it need not necessarily include only quadratic
forms~\cite[]{Donoho+89}.

In multiple Cartesian dimensions, in practice: in the plane, the
space-frequency localization problem has received remarkably little
attention beyond the initial treatment by Slepian
himself~\cite[]{Slepian64}, who restricted his attention to
concentration over circular disks in space and spectral
space~\cite[see also][]{Brander+86,Vandeville+2002}. In this
situation, as in the  one-dimensional case, the concentration operator
commutes with a second-order Sturm--Liouville differential operator,
which greatly facilitates the (numerical) analysis. The scarce recent
work on two-dimensional Slepian functions has focused on their being
amenable to generalized Gaussian
quadrature~\cite[]{DeVilliers+2003,Ma+96,Shkolnisky2007} or on using
them for multiscale out-of-sample extensions~\cite[]{Coifman+2006}
without straying from circular regions of interest. On square or
rectangular domains Slepian functions formed by outer products of
pairs of the pswf~\cite[]{Borcea+2008,Hanssen97} have correspondingly
square or rectangular concentration regions in spectral space, which
is undesirable in geophysical applications~\cite[]{Simons+2003a}.

Geographical regions are not typically squares, circles or rectangles
--- although usually we will require the spectral support of the
Slepian functions to remain disk-like to enable isotropic feature
extraction~\cite[]{Zhang94}. In those studies where Cartesian domains of
arbitrary geometry have been implicit or explicitly considered, as for
applications in image processing~\cite[]{Zhou+84},
radio-astronomy~\cite[]{Jackson+91},
medical~\cite[]{Lindquist+2006,Walter+2008,Yang+2002}, radar or seismic
imaging~\cite[]{Borcea+2008} or in the specific context of
multi-dimensional spectral analysis~\cite[]{Bronez88,Liu+92}, the
Slepian functions are obtained directly from the defining equation,
i.e. by numerical diagonalization of the discretized space-bandwidth
projection operator. This can be implemented
explicitly~\cite[]{Percival+93} or by iterated filtering and
bounding~\cite[]{Jackson+91} as in the
Papoulis~\cite[]{Kennedy+2008,Papoulis75} or Lanczos~\cite[]{Golub+89}
algorithms. Given the characteristic step-shaped eigenvalue spectrum
of the operators~\cite[]{Daubechies88a,Daubechies92,Slepian76,Slepian83}
such procedures are not uniformly stable~\cite[]{Bell+93}, and for
general, non-symmetric domains, no better-behaved commuting operator
exists~\cite[]{Brander+86,Grunbaum+82,Parlett+84}.
In this case we prefer the numerical solution of the Fredholm integral
equation~\cite[]{Tricomi70} with which the concentration problem is
equivalent, by the Nystr\"om method~\cite[]{Nystrom30}, using classical
Gauss-Legendre integration~\cite[]{Press+92}. When both the spatial and
spectral regions of concentration are irregular in shape and of
complete generality, diagonalization of the projection operator is our
only recourse, as we illustrate.

As we have in a sense come full circle as regards ``Slepian's
problem'', at least concerning the scalar case, and because of its
deep connections within the framework of reproducing-kernel Hilbert
spaces~\cite[]{Aronszajn50,Simons2010,Yao67}, our contribution has the
character of a review. We take the reader from the beginnings of the
theory in one linear dimension~\cite[]{Slepian+61} to its generalization
on the surface of the unit sphere~\cite[]{Simons+2006b,Simons+2006a},
and back to the two Cartesian dimensions of the title of this paper,
where, after~\cite{Slepian64}, there remained some work to be
done. Our discussion is only as technical as required to make the
theory available for applications in the natural sciences. Focusing on
the practical and the reproducible, we are distributing all of our
\textsc{Matlab} computer routines freely over the World Wide Web.

\section{Spatiospectral concentration on the real line}\label{realline}

We use~$t$ to denote time or one-dimensional space and~$\omega$ for
angular frequency, and adopt a normalization convention~\cite[]{Mallat98}
in which a real-valued time-domain signal~$f(t)$ and its Fourier
transform~$F(\omega)$ are related by 
\be
\label{slepian0}
f(t)=\jnorm\intinf 
F(\omega)e^{i\omega t}\,d\omega,\qquad 
F(\omega)=\intinf 
f(t)e^{-i\omega t}\,dt.
\ee
Following Slepian, Landau and Pollak~\cite[]{Landau+61,Slepian+61}, the
strictly bandlimited signal   
\be
\label{gband}
g(t)=\jnorm\intW 
G(\omega)e^{i\omega t}\,d\omega
,
\ee
that is maximally concentrated within $|t|\le T$ is the one that
maximizes the ratio
\be
\lambda=\fracd{\int_{-T}^{T}g^2(t)\,dt}
{\intinf g^2(t)\,dt}
.
\label{slepian1}
\ee
Bandlimited functions~$g(t)$ satisfying problem~(\ref{slepian1}) have
spectra~$G(\omega)$ that satisfy 
\begin{subequations}
\label{slepian2}
\be
\intW D(\omega,\omega')\,G(\omega')\,d\omega'
=\lambda\hspf G(\omega),\quad |\omega|\le W,
\ee
\be
D(\omega,\omega')=\frac{\sin [T(\omega-\omega')]}
{\pi (\omega-\omega')}.
\ee
\end{subequations}
The corresponding time- or spatial-domain formulation is 
\begin{subequations}
\label{dude1}
\be
\label{eig1}
\int_{-T}^{T}\Dttp\,g(t')\,dt'=\lambda\hspf g(t),\quad
t\in\mathbb{R}
,
\ee
\be
\label{slepian4}
\Dttp=\frac{\sin [W(t-t')]}{\pi (t-t')}
.
\ee
\end{subequations}
When the prolate spheroidal wave functions $g_1(t),g_2(t),\ldots$ that
solve eq.~(\ref{dude1}) are made orthonormal over $|t|\le\infty$ they
are also orthogonal over the interval $|t|\le T$:   
\be
\intinf g_{\alpha}g_{\beta}\,dt=\dab
,\qquad\int_{-T}^{T}g_{\alpha}g_{\beta}\,dt=
\lambda_{\alpha}\hspf \dab.
\label{slepian8}
\ee
It now follows directly that 
\be\label{nearly1}
\Dttp=\sum_{\alpha=1}^{\infty}g_\alpha(t)g_\alpha(t')\also
\sum_{\alpha=1}^{\infty}g^2_\alpha(t)=\frac{W}{\pi}
.
\ee
A change of variables and a scaling transform eq.~(\ref{slepian2})
into the dimensionless eigenvalue problem
\begin{subequations}
\label{slepian5}
\be
\int_{-1}^{1}D(x,x')\,\psi(x')\,dx'=\lambda\hspf \psi(x),
\ee
\be
D(x,x')=\fracd{\sin [TW(x-x')]}{\pi (x-x')}.
\ee
\end{subequations}
The eigenvalues $\lambda_1>\lambda_2>\ldots$ and eigenfunctions
$\psi_1(x),\psi_2(x),\ldots$ depend only upon the 
time-bandwidth product~$TW$. The sum of the concentration
values~$\lambda$ relates to this as
\be\label{shannon1}
N\D=\sum_{\alpha=1}^{\infty}\lambda_{\alpha}
=\int_{-1}^{1}D(x,x)\,dx
=\frac{2TW}{\pi}
=\frac{(2T)(2W)}{2\pi}. 
\ee
The spectrum of eq.~(\ref{slepian5}) has near-unity and near-zero
eigenvalues separated by a narrow transition
band~\cite[]{Landau65,Slepian+65}. Thus, $N\D$, the ``Shannon number'',
roughly equals the number of significant eigenvalues. In other
words~\cite[]{Landau+62}, it is the approximate dimension of the space
of signals that can be simultaneously well concentrated into finite
time and frequency intervals. The eigenvalue-weighted sum of the
eigenfunctions will be nearly equal to the constant in
eq.~(\ref{nearly1}) over the region of concentration, as
\be
\sum_{\alpha=1}^{\infty}\lambda_{\alpha}\hspf 
g_{\alpha}^2(t)\approx
\sum_{\alpha=1}^{N\D}\lambda_{\alpha}\hspf 
g_{\alpha}^2(t)\approx
\left\{\barray{ll}
N\D/(2T) & \mbox{if $-T\le t \le T$,}\\
0 & \mbox{otherwise.}
\earray\right.
\label{sumofsq1}
\ee

The integral operator in eq.~(\ref{slepian5}) commutes with a
Sturm--Liouville differential operator~\cite[]{Slepian83,Slepian+61}, to
the effect that although~$\chi\ne\lambda$, we can solve for the
functions~$\psi$ also from 
\be
\label{sturmliouville1}
\frac{d}{dx}\left[(1-x^2)\frac{d\psi}{dx}\right]+
\left[\chi-\frac{(N\D)^2\pi^2}{4}x^2\right]  
\hsp\psi=0,\qquad |x|\le 1.
\ee
At~$N$ discrete values of~$x=0,...,N-1$, the eigenfunctions~$\psi(x)$ of
eq.~(\ref{sturmliouville1}) can be found by diagonalization of a simple
symmetric tridiagonal matrix~\cite[]{Grunbaum81a,Percival+93,Slepian78}
with elements 
\begin{subequations}\label{trione}\ber
T_{xx}&=&[(N-1-2x)/2]^2\cos(2\pi W),\\
T_{x\,x+1}&=&(x+1)(N-x-1)/2
.
\eer
\end{subequations}
The matching concentration eigenvalues~$\lambda$ can then be obtained
directly from eq.~(\ref{slepian5}). Compared to those, the
Sturm--Liouville eigenvalues~$\chi$ are very regularly distributed and
thus the computation of the Slepian functions via diagonalization of
eq.~(\ref{trione}) is always stable~\cite[]{Percival+93}. 

\section{Spatiospectral concentration on a sphere}

We use~$\rhat=(\theta,\phi)$ to denote a location at
colatitude~$\theta$ and longitude~$\phi$ on the unit
sphere~$\Omega=\{\rhat: \|\rhat\|=1\}$, and adopt a normalization
convention~\cite[]{Dahlen+98,Edmonds96} in which a  real-valued
time-domain signal~$f(\rhat)$ and its spherical-harmonic
transform~$f_{lm}$ at degree~$l$ and order~$m$  are related by 
\be
\label{expansion}
f(\rhat)=\sumsh {f}_{lm}\Ylmrh,\qquad {f}_{lm}=\into
f(\rhat)\hspace*{0.05em}\Ylmrh\domg
.
\ee
Following Simons, Wieczorek and Dahlen~\cite[]{Simons+2006a}
the strictly bandlimited signal   
\be\label{bandlg}
\grh=\sumshL \glm \Ylmrh
\ee
that is maximally concentrated within a region~$R\subset\Omega$ 
is the one that maximizes the ratio
\be
\lambda=\fracd{\intr g^2(\rhat)\domg}{\into^{}g^2(\rhat)\domg}
.
\label{normratio}
\ee
Bandlimited functions~$\grh$ satisfying problem~(\ref{normratio}) have
spectra~$\glm$ that satisfy
\begin{subequations}
\label{fulleigen1}
\be
\label{halfeigen}
\sumshLp\Dlmlmp
\glmp=\lambda\hspf \glm ,\qquad 0\le l\le L,
\ee
\be
\Dlmlmp=\intr\Ylm\Ylmp\domg.
\label{Dlmlmpdef}
\ee
\end{subequations}
Through the addition theorem~\cite[]{Dahlen+98}, the corresponding
spatial-domain formulation is 
\begin{subequations}
\label{firsttimeint}
\be
\label{eig3}
\intr  \Drhrhp \,\grhp\domg'=
\lambda\hspf \grh,\qquad \rhat\in\Omega,
\ee
\be
\Drhrhp
=\sum_{l=0}^L\tlofp\!P_l(\rhat\cdot\rhat'),
\label{banddelta}
\ee
\end{subequations}
where~$P_l$ is the Legendre function. When the 
functions
$g_1(\rhat),g_2(\rhat), \ldots,g_{(L+1)^2}(\rhat)$  that
solve eq.~(\ref{firsttimeint}) are orthonormal over $\Omega$
they are also orthogonal over the region~$R$:
\be
\into g_{\alpha}g_{\beta}\domg=\dab,\qquad
\intr  g_{\alpha}g_{\beta}\domg=\lambda_{\alpha}\dab.
\label{orthog}
\ee
From eqs.~(\ref{firsttimeint})--(\ref{orthog}) we then immediately obtain
the relations
\be
\Drhrhp=\sum_{\alpha=1}^{\Lpot}g_{\alpha}(\rhat)g_{\alpha}(\rhat')\also
\label{newbasis}
\sum_{\alpha=1}^{\Lpot}g_{\alpha}^2(\rhat)
=\frac{\Lpot}{4\pi}.
\label{sumofsq}
\ee

Asymptotically, as the spatial area 
$A\rightarrow 0$, and~$L\rightarrow\infty$, eq.~(\ref{firsttimeint}) 
becomes~\cite[]{Simons+2007,Simons+2006a}
\begin{subequations}
\label{3Dscaled}
\be
\int_{R_{*}}\!\!\Dxixip
\,\psi(\bxi')\,d\Omega'_*=\lambda\hsp\psi(\bxi),
\ee
\be
\label{asymptotic}
\Dxixip=\frac{(L+1)\sqrt{A/4\pi}}{2\pi}
\fracd{J_1[(L+1)\sqrt{A/4\pi}\,\,\|\bxi-\bxi'\|]} 
{\|\bxi-\bxi'\|},
\ee
\end{subequations}
where the scaled region~$R_{*}$  now has area~$4\pi$ and~$J_1$ 
is the first-order Bessel function of the first kind. As in the
one-dimensional case~(\ref{shannon1}), the eigenvalues 
$\lambda_1\ge\lambda_2\ge\ldots$ and eigenfunctions 
$\psi_1(\bxi),\psi_2(\bxi),\ldots$ depend only upon the product of the
maximal degree~$L$ and the area~$A$.  The sum of the concentration
values~$\lambda$ , the space-bandwidth product or ``spherical Shannon
number'', $N\DDD$, once again roughly the number of
significant eigenvalues, is: 
\ber
\label{shannon2}
N\DDD&=&\sum_{\alpha =1}^{\Lpot}\lambda_{\alpha}=
\sum_{l=0}^L\sum_{m=-l}^l 
D_{lm,lm}=\int_RD(\rhat,\rhat)\,d\Omega,\\
&=&\int_{R_{*}}\!\!\Dxixi
\,d\Omega_*
=\Lpot\,\frac{A}{4\pi}
\nnr
.
\eer
The first~$N\DDD$ orthogonal eigenfunctions
$g_{\alpha},\alpha=1,2,\ldots,N\DDD$, with significant eigenvalues
$\lambda_{\alpha}\approx 1$, provide an essentially uniform coverage
of the region~$R$, reaching the constant of eq.~(\ref{sumofsq}) in the
following sense:
\be
\sum_{\alpha=1}^{\Lpot}\lambda_{\alpha}\hspf 
g_{\alpha}^2(\rhat)\approx
\sum_{\alpha=1}^{N\DDD}\lambda_{\alpha}\hspf 
g_{\alpha}^2(\rhat)\approx
\left\{\barray{ll}
N\DDD/A & \mbox{if $\rhat\in R$},\\
0 & \mbox{otherwise.}
\earray\right.
\label{sumofsq2}
\ee
Irrespectively of the particular region of concentration that they
were designed for, the complete set of bandlimited spatial Slepian
eigenfunctions $g_1,g_2, \ldots,g_{(L+1)^2}$ is a basis for
bandlimited scalar processes anywhere on the surface of the unit
sphere~\cite[]{Simons+2006b,Simons+2006a}.  The reduced set
$g_{1},\ldots, g_{N\DDD}$, with eigenvalues $\lambda\approx 1$, is 
an approximate basis for bandlimited processes that are primarily
localized to the region~$R$~\cite[]{Simons+2009b}. Thus $N\DDD$ is the
approximate dimension of the space of signals that can be
simultaneously well concentrated into finite spatial and spectral
intervals on the surface of the unit sphere. Eq.~(\ref{3Dscaled})
depends only on this combination of bandwidth and spatial area, as in
the one-dimensional case, eq.~(\ref{slepian5}).  

When the region of concentration is a circularly symmetric cap of
colatitudinal radius~$\Theta$, centered on the North Pole, the
colatitudinal parts~$\sfg(\theta)$ of the separable solutions to
eq.~(\ref{firsttimeint}), 
\be
g(\theta,\phi)=\left\{
\barray{l@{\quad\mbox{if}\hspace{0.6em}}l}
\rule[-2mm]{0mm}{6mm}\sqrt{2}\,\sfg(\theta)\cos (m\phi) & -L\le m<0,\\
\rule[-2mm]{0mm}{6mm}\sfg(\theta)                     & m=0,\\
\rule[-2mm]{0mm}{6mm}\sqrt{2}\,\sfg(\theta)\sin (m\phi) & 0< m\le L,\\
\earray
\right.
\label{polarg2}
\ee
are, owing to a commutation relation~\cite[]{Grunbaum+82}, identical to
those of the Sturm--Liouville equation
\ber\label{sturmliouville2}
\lefteqn{\fracd{d}{d\mu}\left[(\mu-\cos\Theta)
(1-\mu^2)\fracd{d\sfg}{d\mu}\right]}\hspace{6em}
\\
\nnr
&&{}+\left[\chi+
L(L+2)\mu
-\fracd{m^2(\mu-\cos\Theta)}{1-\mu^2}
\right]\sfg=0,
\eer
with~$\mu=\cos\theta$ and~$\chi\ne\lambda$. At constant~$m$, their
values $\glm$  in the expansion~(\ref{bandlg}) can be found
by diagonalization of a simple symmetric tridiagonal
matrix~\cite[]{Grunbaum+82,Simons+2006a} with elements
\begin{subequations}
\label{gdefi1}
\ber 
T_{ll}&=&-l(l+1)\cos\Theta,\\
T_{l\,l+1}&=&\big[l(l+2)-L(L+2)\big]
\sqrt{\fracd{(l+1)^2-m^2}{(2l+1)(2l+3)}}.
\eer
\end{subequations}
When the region of concentration is a pair of axisymmetric
polar caps of common radius~$\Theta$ centered on the
North and South Pole, the~$\sfg(\theta)$ solve the Sturm--Liouville
equation 
\ber\label{sturmliouville22}
\lefteqn{\fracd{d}{d\mu}\left[(\mu^2-\cos^2\Theta)
(1-\mu^2)\fracd{d\sfg}{d\mu}\right]}\hspace{6em}
\\
\nnr
&&{}+\left[
\chi + L_p(L_p+3)\mu^2
-\fracd{m^2(\mu^2-\cos^2\Theta)}{1-\mu^2}
\right]\sfg=0,
\eer
where~$L_p=L$ or~$L_p=L-1$ depending whether the order~$m$ of the
functions~$g(\theta,\phi)$ in eq.~(\ref{polarg2}) is odd or
even, and whether the bandwidth~$L$ is odd or
even~\cite[]{Grunbaum+82,Simons+2006b}. The expansion coefficients of
the optimally concentrated antipodal polar-cap eigenfunctions require the
numerical diagonalization of the symmetric tridiagonal
matrix~\cite[]{Simons+2006b} with elements
\begin{subequations}
\label{gdefi}
\ber
T^p_{ll}&=&-l(l+1)\cos^2\Theta
+\frac{2}{2l+3}\left[(l+1)^2-m^2\right]\\
 \nnr&&{}+[(l-2)(l+1)-L_p(L_p+3)]
\left[\frac{1}{3}-\frac{2}{3}\,
\fracd{3m^2-l(l+1)}{(2l+3)(2l-1)}\right],\\ 
T^p_{l\,l+2}&=&\fracd{\big[l(l+3)-L_p(L_p+3)\big]}{2l+3}
\sqrt{\fracd{\left[(l+2)^2-m^2\right]
\left[(l+1)^2-m^2\right]}{(2l+5)(2l+1)}}.\label{skip}
\eer
\end{subequations}
Every other degree in the expansion for the equatorially
(anti-)symmetric double-polar cap functions is
skipped~\cite[]{Simons+2006b}, hence the $l+2$ subscript for the
elements off the main diagonal in eq.~(\ref{skip}). The 
concentration values~$\lambda$ can be determined from the defining
equations~(\ref{fulleigen1}) or~(\ref{firsttimeint}). Compared to
these, the Sturm--Liouville eigenvalues~$\chi$ in
eqs.~(\ref{sturmliouville2}) and~(\ref{sturmliouville22}) are very
regularly spaced, and thus the computations for these special cases
are inherently stable.

\section{Spatiospectral concentration in the Cartesian plane}
\label{cartesian}

We now turn to the multidimensional Cartesian case, first discussed
by~\cite{Slepian64}, noting that we have set ourselves up for 
a result that is analogous to eq.~(\ref{3Dscaled}). Indeed, in the
asymptotic regime of an infinitely large bandwidth and an
infinitesimally small region on the surface of the unit sphere, the
spherical and Cartesian concentration problems are of course
equivalent~\cite[]{Simons+2007,Simons+2006a}. 

Another note concerns the equivalence between the temporal or spatial
and spectral forms of the concentration problems. In writing
eqs.~(\ref{slepian2})--(\ref{dude1}), and
eqs.~(\ref{fulleigen1})--(\ref{firsttimeint}), respectively, we have
exclusively considered strictly bandlimited, time- or
space-concentrated Slepian functions. Strictly time- or spacelimited,
band-concentrated functions can be obtained via an appropriate
restriction of the integration domains in these equations, and the
resulting new functions can be obtained from the old ones by simple
truncation and rescaling. Both for the one-dimensional and spherical
situations this distinction is usually made more explicitly
elsewhere~\cite[]{Landau+61,Simons+2006a,Slepian+61}, and in what
follows we once again treat both cases separately.

For applications in the geosciences, we focus on the two-dimensional,
flat geometry of geographical maps: the Cartesian plane
$\mathbb{R}^2$, with a spatial coordinate vector~$\bx$ and a spectral
coordinate vector~$\bk$. To wit, we formulate the concentration
problem as follows: in spite of the Paley-Wiener theorem, which states
that functions cannot be spatially and spectrally restricted at the
same time~\cite[]{Daubechies92,Mallat98}, can we construct (real)
functions that are localized to, say, the shape of Belgium, in
$\bx$-space, while having a Fourier transform localized to, say, the
(Hermitian) shape of a pair of triangles, in $\bk$-space? Yes we can!  

\subsection{Preliminary considerations}

A real-valued, square-integrable function~$f(\bx)$ defined in
the plane has the two-dimensional Fourier representation
\be
\label{Bfourier}
f(\bx)=\fnorm\intinft
F(\bk)e^{i\bk\cdot\bx}\dbk,\qquad
F(\bk)=\intinft
f(\bx)e^{-i\bk\cdot\bx}\dbx,
\ee
with~$F(\bk)=F^*(-\bk)$. The Fourier orthonormality relation is
\be
\label{BFortho}
\fnorm\intinft
e^{i\bk\cdot(\bx-\bx')}\dbk
=\delta(\bx,\bx')
,
\ee
which defines the delta function in the usual distributional sense
\be
\intinft f(\bx')\delta(\bx,\bx')\dbx'=f(\bx)
.
\ee
Likewise, in the spectral domain we may write
\be
\label{BFortho2}
\fnorm\intinft
e^{i(\bk-\bk')\cdot\bx}\dbx
=\delta(\bk,\bk')
.
\ee
By Parseval's relation the energies in the spatial and spectral
domains are identical:
\be
\label{Btwopowers}
\intinft f^2(\bx)\dbx=\fnorm\!\!
\intinft |F(\bk)|^2\dbk.
\ee

\subsection{Spatially concentrated bandlimited functions}
\label{theg}

We use~$g(\bx)$ to denote a real function that is bandlimited
to~$\ssK$, an arbitrary subregion of spectral space, 
\be
\label{Bgdefn}
g(\bx)=\fnorm\!\!
\intK G(\bk)e^{i\mbf{k}\cdot\bx}\dbk,
\ee
Following~\cite{Slepian64}, we seek to concentrate the power
of~$g(\bx)$ into a finite spatial region~$\ssR\subset\mathbb{R}^2$,
of area~$A$, by maximizing the energy ratio
\be
\label{Bnormratio}
\lambda=\fracd{\intR g^2(\bx)\dbx}
{\intinft g^2(\bx)\dbx}.
\ee
Upon inserting the representation~(\ref{Bgdefn}) of~$g(\bx)$ into
eq.~(\ref{Bnormratio}) we can express the concentration in the form of
the Rayleigh quotient 
\be
\label{Bnormratio2}
\lambda=\fracd{\intK\intK G^*(\bk)
\Dkkp G(\bk')\dbk
\dbk'}
{\intK |G(\bk)|^2\dbk}
,
\ee
where we have used Parseval's relation~(\ref{Btwopowers}) and defined
the positive-definite quantity 
\be
\label{Beigen2}
\Dkkp=\fnorm\!\!
\intR e^{i(\bk'-\bk)\cdot\bx}\dbx,
\ee
which is Hermitian, $\Dkkp=\Dkpkc$. Bandlimited functions~$g(\bx)$
that maximize eqs.~(\ref{Bnormratio})--(\ref{Bnormratio2}) solve the 
Fourier-domain Fredholm integral equation 
\be
\label{Beigen1}\label{Beigen12}
\intK \Dkkp\,G(\bk')\dbk'
=\lambda\hsp G(\bk),\qquad\bk\in\ssK
.
\ee
Comparison of eq.~(\ref{Beigen2}) with eq.~(\ref{BFortho2}) leads to
the interpretation of the spectral-domain kernel~$\Dkkp$ as a
spacelimited spectral delta function. We rank order the concentration
eigenvalues so that $1>\lambda_1\ge\lambda_2\ge\ldots>0$. Upon
multiplication of eq.~(\ref{Beigen1}) with~$e^{i\bk\cdot\bx}$  
and integrating over all $\bk\in\ssK$, we obtain
the corresponding  problem in the spatial domain as 
\begin{subequations}
\label{Beigen34}\label{eig2}
\be
\label{Beigen3}
\intR\! \Dxxp\,g(\bx')\dbx'
=\lambda\hspo g(\bx),\quad\bx\in\mathbb{R}^2,
\ee
\be
\label{Beigen4}
\Dxxp=\fnorm
\intK e^{i\bk\cdot(\bx-\bx')}\dbk.
\ee
\end{subequations}
Comparison of eq.~(\ref{Beigen4}) with eq.~(\ref{BFortho}) shows that
the Hermitian spatial-domain kernel $\Dxxp=D^*(\bx',\bx)$ is a
bandlimited spatial delta function.  The bandlimited spatial-domain
eigenfunctions $g_1(\bx),g_2(\bx),\ldots$ may be chosen to be
orthonormal over the whole plane~$\mathbb{R}^2$, in which case they are
also orthogonal over the region~$\ssR$:  
\be
\intinft
g_{\alpha}g_{\beta}\dbx=\dab,\qquad
\intR
g_{\alpha}g_{\beta}\dbx=\lambda_{\alpha}\dab.
\label{Borthog1}
\label{Borthog2}
\ee
In this normalization the eigenfunctions of eq.~(\ref{Beigen34})
represent the kernel in~(\ref{Beigen4}) as  
\be
\label{Mercer}
\Dxxp=\sum_{\alpha=1}^\infty g_\alpha(\bx) g_\alpha(\bx')
.
\ee
This form of Mercer's theorem~\cite[]{Flandrin98,Tricomi70}
is verified by substituting the right hand side of
eq.~(\ref{Mercer}) into eq.~(\ref{Beigen3}) and using the
orthogonality~(\ref{Borthog2}). 

\subsection{Spectrally concentrated spacelimited functions}
\label{theh}

We use~$h(\bx)$ to denote a function that is spacelimited,
i.e. vanishes outside the arbitrary region~$\ssR$ of
physical space:
\be\label{Khdefn}
H(\bk)=\intR
h(\bx)e^{-i\bk\cdot\bx}\dbx. 
\ee
To concentrate the energy of~$h(\bx)$ into the finite
spectral region~$\ssK$, we maximize
\be
\label{Knormratio}
\lambda=\fracd{\intK |H(\bk)|^2\dbk}
{\intinft |H(\bk)|^2\dbk}.
\ee
Upon using eq.~(\ref{Btwopowers}) and inserting the
representation~(\ref{Khdefn}) of $H(\bk)$ into eq.~(\ref{Knormratio})
we can rewrite $\lambda$ in the form
\be
\label{Knormratio2}
\lambda=\fracd{\intR\intR h(\bx)
\Dxxp h(\bx')\dbx
\dbx'}
{\intR h^2(\bx)\dbx},
\ee
where we again encounter the quantity~(\ref{Beigen4}),
\be
\label{Keigen}
\Dxxp=\fnorm
\intK e^{i\bk\cdot(\bx-\bx')}\dbk
.
\ee
Once again by Rayleigh's principle, spacelimited functions~$h(\bx)$
that maximize the quotient~$\lambda$ in
eqs.~(\ref{Knormratio})--(\ref{Knormratio2}) solve the 
spatial-domain Fredholm integral equation 
\be
\label{Keigen1}
\intR \!\Dxxp\,h(\bx')\dbx'
=\lambda\hspo h(\bx),\quad \bx \in \ssR
.
\ee
Eq.~(\ref{Keigen1}) is identical to eq.~(\ref{Beigen34}) save for the 
restriction to the domain~$\ssR$. The eigenfunctions~$h(\bx)$ that
maximize the spectral norm ratio~(\ref{Knormratio}) are identical,
within the region~$\ssR$, to the eigenfunctions~$g(\bx)$ that
maximize the spatial norm ratio~(\ref{Bnormratio}). The associated
eigenvalues $1>\lambda_1\ge\lambda_2\ge\ldots>0$  are a measure both
of the spatial concentration of~$g(\bx)$ within the region~$\ssR$ and
of the spectral concentration of~$h(\bx)$ to the wave
vectors~$\bk\in\ssK$. Identifying 
\be
h(\bx)=\left\{\barray{ll}
g(\bx) & \mbox{if $\bx\in \ssR$},\\
0 & \mbox{otherwise},
\earray\right.
\label{hequalsg}
\ee
the normalization is such that
\be
\intinft
h_{\alpha}h_{\beta}\dbx=\intR h_{\alpha}h_{\beta}\dbx=
\lambda_\alpha\dab.
\label{Borthoh1}
\ee
The null-space consisting of all non-bandlimited functions that
have both no energy outside~$\ssR$ nor inside~$\ssK$ is
of little consequence to us.

\section{Slepian Symmetry}\label{slepsym} 

Under what has been called ``Slepian symmetry''~\cite[]{Brander+86} the
solutions to eqs.~(\ref{Beigen34}) and~(\ref{Keigen1}) take on
particularly attractive analytic forms. Anticipating a switch to
polar coordinates $\bx=(r,\theta)$ we introduce~$J_m(k)$, the Bessel
function of the first kind and of integer order~$m$
\cite[]{Abramowitz+65,Gradshteyn+2000}. The Bessel functions satisfy the
symmetry condition $J_{-m}(k)=(-1)^mJ_m(k)$, the relation 
\be
1=J_0^2(k)+2\suml_{m=1}^{\infty}J_m^2(k)
\label{additionBJ}
\ee
and the identity~\cite[]{Jeffreys+88} 
\be
J_0(k\|\bx-\bx'\|)=J_0(kr)J_0(kr')
+2\suml_{m=1}^{\infty}J_m(kr)J_m(kr')\cos [m(\theta-\theta')]
\label{jeffs}
.
\ee
Furthermore, we have the particular formulas
\ber
\label{jzero}\label{jone}
\lefteqn{\hspace{-18em}J_0(k)=\jnorm\intth e^{ik\cos\theta}\dth,\qquad
J_1(k)=k^{-1}\int_0^{k} J_0(k')\,k'\dk',}\\
\label{jhalf}
J_{1/2}(k)=\sqrt{\frac{2}{\pi k}}\,\sin k,
\eer
the derivative identity 
$\frac{d}{dk}[kJ_1(k)]=kJ_0(k)$
and the limits
\be
\label{jlim}
\lim_{k\rar 0} \frac{J_1(k)}{k}\rar \half
\also
\lim_{k\rar 0} \frac{J_m(k)}{\sqrt{k}}\rar 0
.
\ee

\subsection{Circular bandlimitation}
\label{circband}

Limitation to the disk-shaped $\ssK=\{\bk:\|\bk\|\le K\}$ allows us
to rewrite the spatial kernel of eq.~(\ref{Keigen}) using polar coordinates as
follows. Letting $k=\|\bk\|=(\bk\cdot\bk)^{1/2}$ and~$\theta$ be the
angle between the wave vector~$\bk$ and the vector $\bx-\bx'$, 
\be
\bk\cdot(\bx-\bx')=k\,\|\bx-\bx'\|\cos\theta,
\ee
and using eq.~(\ref{jzero}), we obtain an expression alternative to
eq.~(\ref{Keigen}) as   
\begin{subequations}\ber
\label{Bbesskern}
\Dxxp&=&\fnorm\intk\intth e^{ik\|\bx-\bx'\|\cos\theta}\dth\,k\dk,
\\
&=&\jnorm\intk\!J_0(k\|\bx-\bx'\|)\,k\dk,
\label{Bbesskern1}\\
&=&\fracd{K\hspo 
J_1(K\|\bx-\bx'\|)}{2\pi\, \|\bx-\bx'\|}.\label{Bbesskern2}
\eer
\end{subequations}
We notice the symmetry $\Dxxp=D(\bx',\bx)=D(\|\bx-\bx'\|)$ and the
equivalence of eq.~(\ref{Bbesskern2}) with the asymptotic
expression~(\ref{asymptotic}), as
expected~\cite[]{Simons+2007,Simons+2006a}. With~\cite{Coifman+2006}
and through eq.~(\ref{jhalf}), we furthermore note that both
eqs.~(\ref{Bbesskern2}) and~(\ref{slepian4}) are $n$-dimensional
versions of the general form  
$
[K/(2\pi)]^{n/2} J_{n/2}(K\|\bx-\bx'\|)/(\|\bx-\bx'\|^{n/2})$,
for~$n=2$, and $n=1$ (where $K=W)$, respectively.
As eq.~(\ref{jlim}) now shows,
\be\label{Dxx}
\Dxx=\frac{K^2}{4\pi}
,
\ee
and thus eqs.~(\ref{Borthog2})--(\ref{Mercer}) allow us to write the
sum of the concentration values $\lambda$ as the ``planar Shannon  number''
\be
\label{shannon3}
N\DD=\sum_{\alpha=1}^{\infty}\lambda_{\alpha}=\intR \!\Dxx\dbx
=K^2\,\fracd{A}{4\pi}
=\fracd{(\pi K^2)(A)}{(2\pi)^2}
,
\ee
where~$A$ is the area of the spatial region of concentration
$\ssR$. As in the one-dimensional case, eq.~(\ref{shannon1}),
the Shannon number is equal to the area of the function in the
spatiospectral phase plane times the ``Nyquist density'',
as expected~\cite[]{Daubechies88a,Daubechies92}. Roughly
equal to the number of significant eigenvalues, $N\DD$ again is the
effective dimension of the space of ``essentially'' space- and bandlimited
functions in which the reduced set of two-dimensional functions
$g_1,g_2,\dots,g_{N\DD}$ may act as an efficient orthogonal basis. 
Using eqs.~(\ref{Dxx}) with~(\ref{Mercer}) we furthermore see that the
sum of the squares of all of the bandlimited eigenfunctions,
independent of position $\bx$ in the plane, is the constant 
\be
\sum_{\alpha=1}^\infty g^2_\alpha(\bx)=\frac{K^2}{4\pi}
.
\ee
Likewise, since the first~$N\DD$ eigenfunctions $g_1,g_2,\ldots,g_{N\DD}$
have eigenvalues near unity and lie mostly within $\ssR$, and
the remainder $g_{N\DD+1},g_{N\DD+2},\ldots,g_{\infty}$ have eigenvalues
near zero and lie mostly outside $\ssR$ in $\mathbb{R}^2$,
we expect the eigenvalue-weighted sum of squares to be  
\be\label{heuristic}
\sum_{\alpha=1}^\infty \lambda_\alpha g^2_\alpha(\bx)
\approx
\sum_{\alpha=1}^{N\DD} \lambda_\alpha g^2_\alpha(\bx)
\approx
\label{sumofsq3}\left
\{\barray{ll}
N\DD/A & \mbox{if $\bx\in \ssR$},\\
0 & \mbox{otherwise.}
\earray\right.
\ee
This heuristic finding is of great importance in the analysis and
representation of signals using the Slepian functions as an approximate
basis, as much as for their use as tapers to perform spectral analysis
on data from which we thus expect to extract all relevant statistical
information with minimal loss or leakage near the edges of the region
under consideration~\cite[]{Walden90a}.  

\ssec{Scaling analysis I}

Introducing scaled independent and dependent variables
\be
\label{Bhuge}
      \bxi=\sqrt{\fracd{4\pi}{A}}\,\bx,\qquad
\bxi'=\sqrt{\fracd{4\pi}{A}}\,\bx',
\qquad
\psi(\bxi)=g(\bx),\qquad
\psi(\bxi')=g(\bx'),
\ee
we can rewrite eqs.~(\ref{Beigen34}) and~(\ref{Bbesskern2}) as
\begin{subequations}
\label{BDscaled}
\be
\int_{\ssR_{*}}\!\!\Dxixip
\,\psi(\bxi')\,d\bxi'=\lambda\hspo \psi(\bxi),
\label{Bscaled}
\ee
with~$\ssR_{*}$, of area~$4\pi$, the image of the concentration
region~$\ssR$ under~(\ref{Bhuge}), and 
\be\label{Dscaled}
\Dxixip=\frac{\sqrt{N\DD}}{2\pi}
\fracd{J_1(\sqrt{N\DD}\|\bxi-\bxi'\|)} 
{\|\bxi-\bxi'\|}.
\ee
\end{subequations}
Eq.~(\ref{BDscaled}) reveals that the eigenvalues
$\lambda_1,\lambda_2,\ldots$ and the scaled eigenfunctions
$\psi_1(\bxi)$, $\psi_2(\bxi),\ldots$ depend on the maximum circular
bandwidth~$K$ and the spatial concentration area~$A$ only through the
planar Shannon number~$N\DD$, as they do in the one-dimensional,
eq.~(\ref{slepian5}), and asymptotic spherical, eq.~(\ref{3Dscaled}),
cases. Under this scaling the Shannon number remains 
\be
N\DD=\int_{\ssR_{*}}\!
\Dxixi\,d\bxi=\frac{N\DD}{4\pi}\int_{\ssR_{*}}\!\!\,d\bxi=N\DD
.
\ee 

\subsection{Circular spacelimitation}

If in addition to the circular spectral limitation, physical space is
also circularly limited, in other words, if the spatial region of
concentration or limitation~$\ssR$ is a circle of radius~$R$,
then a polar coordinate, $\bx=(r,\theta)$, representation 
\be g(r,\theta)=\left\{
\begin{array}{l@{\quad\mbox{if}\hspace{0.6em}}l}
\rule[-2mm]{0mm}{6mm}\sqrt{2}\,\sfg(r)\cos (m\theta) & m<0,\\
\rule[-2mm]{0mm}{6mm}\sfg(r)                     & m=0,\\
\rule[-2mm]{0mm}{6mm}\sqrt{2}\,\sfg(r)\sin (m\theta) & m>0,\\
\end{array}
\right.
\label{Bpolarg}
\ee
may be used to decompose eq.~(\ref{BDscaled}) into an infinite series
of non-degenerate fixed-order eigenvalue problems, one for each
order~$\pm m$.  We do this most easily by inserting eq.~(\ref{jeffs}) into
eq.~(\ref{Bbesskern1}), the obtained result with eq.~(\ref{Bpolarg})
into eq.~(\ref{Beigen3}), and using the orthogonality of the functions
$\ldots,\sqrt{2}\cos (m\theta),\ldots, 1,\ldots, \sqrt{2}\sin
(m\theta),\ldots$ over the polar angles $0\le\theta<2\pi$. The integral
equation for the fixed-$m$ radial function~$\sfg(r)$ is then given by 
\begin{subequations}
\label{Bfixed}\be
\int_{0}^{R}D(r,r')\,\sfg(r')\,r'dr'
={\lambda}\hspo \sfg(r),
\label{BFIXEDM}
\ee
with the fixed-$m$ symmetric kernel
\be
D(r,r')=K^2\!\int_{0}^{1}
J_m(Kp\hspf r)\,
J_m(Kp\hspf r')\,p\hspo dp.
\label{BFIXEDM2}
\ee
\end{subequations}
We rank order the distinct but pairwise occurring eigenvalues obtained
by solving each of the fixed-order eigenvalue
problems~(\ref{Bfixed}) so that $1>\lambda_1>\lambda_2>\cdots >0$,
and we orthonormalize the associated eigenvectors $\sfg_1,\sfg_2,\dots$ as
in eq.~(\ref{Borthog1}) so that
\be
2\pi\intzinf
\sfg_{\alpha}(r)\sfg_{\beta}(r)\,rdr=\dab,
\label{Borthog3}
\qquad 2\pi \int_0^R \sfg_{\alpha}(r)\sfg_{\beta}(r)\,rdr=
\lambda_{\alpha}\dab.
\label{Borthog4}
\ee
We shall denote the radial part of the functions~$h(r,\theta)$ by
$\sfh(r)$ in Sect.~\ref{numex}.

\ssec{Scaling analysis II}

Finally, the scaling transformations
\be
\label{Bmapping2}
\xi=r/R,\qquad\quad\,\,
\xi'=r'\hspace{-0.2em}/R,
\qquad 
\psi(\xi)=\sfg(r),\qquad
\psi(\xi')=\sfg(r'),
\ee
convert eq.~(\ref{Bfixed}) into the scaled eigenvalue problem
\begin{subequations}
\label{BDscaled2}
\be
\label{Bscaled2}
\int_{0}^{1}D(\xi,\xi')\hspo \psi(\xi')\,
\xi'\hspf d\xi'=\lambda\hspo \psi(\xi),
\ee
with the fixed-$m$ kernel
\be
\label{Dscaled2}
D(\xi,\xi')=4N\DD\int_{0}^{1}\!
J_m\big(2\sqrt{N\DD}\,p\hspf \xi\big)\,
J_m\big(2\sqrt{N\DD}\,p\hspf \xi'\big)\,p\hspo dp,
\ee
\end{subequations}
which is dependent only upon the Shannon number 
\be
N\DD=K^2\frac{R^2}{4}
\label{scaleN}
.
\ee
The number of significant eigenvalues per angular order~$m$
is~\cite[]{Simons+2007,Simons+2006a} 
\ber
N\DD_m&=&\int_0^1D(\xi,\xi)\,\xi \hspo d\xi\label{SEEPROOF}
=4N\DD\int_0^1\!\!\int_0^1\!
J_m^2\big(2\sqrt{N\DD}\,p\hspo \xi\big)\,p\hspo dp\,\xi
\hspo d\xi,\label{willuse}\\ 
&=&{}2N\DD\left[J^2_m\big(2\sqrt{N\DD}\,\big)+
J^2_{m+1}\big(2\sqrt{N\DD}\,\big)\right]\label{nsubm3}
\\
&&{}-(2m+1)\sqrt{N\DD}J_m\big(2\sqrt{N\DD}\,\big)
J_{m+1}\big(2\sqrt{N\DD}\,\big)\nnr\\
&&{}-\frac{m}{2}\left[
1-J_0^2\big(2\sqrt{N\DD}\,\big)-2\suml_{n=1}^{m}
J_n^2\big(2\sqrt{N\DD}\,\big)
\right].\nnr
\eer
The complete Shannon number~$N\DD$ is preserved,
inasmuch~as, per eqs.~(\ref{willuse}) and~(\ref{additionBJ}),
\be\label{Ntwice}
N\DD= N\DD_0+2\suml_{m=1}^{\infty}N\DD_m=
4N\DD\int_{0}^{1}\!\!\int_{0}^{1}\,p\hspo dp\,\xi\,d\xi=N\DD
.\ee

\ssec{Sturm--Liouville character and tridiagonal matrix formulation}


As noted by~\cite{Slepian64}, eq.~(\ref{BDscaled2}) is an
iterated version of the equivalent ``square-root'' equation   
\be
2\sqrt{N\DD}\int_0^1\!J_m\big(2\sqrt{N\DD}\,\xi\hspf \xi'\big)\,
\psi(\xi')\,\xi'\hspf d\xi'=\sqrt{\lambda}\hspo \psi(\xi).
\label{slepsqrt}
\ee
To see this, it suffices to substitute for $\psi(\xi')$ on the left
hand side of eq.~(\ref{slepsqrt}) using eq.~(\ref{slepsqrt})
itself. Thus, the eigenvalues $\lambda_1,\lambda_2,\ldots$ and eigenfunctions
$\psi_1(\xi),\psi_2(\xi),\ldots$ of eq.~(\ref{BDscaled2}) may
alternatively be found by solving the equivalent
equation~(\ref{slepsqrt}). A further reduction can be obtained by
substituting a scaled Shannon number or bandwidth and rescaling,
\be\label{rescaling}
c=2\sqrt{N\DD},\qquad \gamma=\sqrt{\lambda}/\sqrt{c},
\qquad  \varphi(\xi)=\sqrt{\xi}\,\psi(\xi),
\ee
to yield the more symmetric form
\be\label{GL}
\int_{0}^{1} J_m (c\hspf \xi\xi') \sqrt{c\hspf \xi\xi'}
\,\varphi(\xi')\,d\xi'
=
\gamma\hspo \varphi(\xi).
\label{1Dcircular}
\ee
On the domain~$0\le\xi \le 1$ the~$\varphi(\xi)$
also solve a Sturm--Liouville equation, 
\be\label{sturmliouville3}
\fracd{d}{d\xi}\left[(1-\xi^2)\frac{d\varphi}{d\xi}\right] + \left( \chi+
\fracd{1/4-m^2}{\xi^2}-4N\DD \xi^2  \right)\varphi =0,\qquad |\xi|\le 1
,
\ee
for some~$\chi\ne\lambda$. When~$m=\pm1/2$ eq.~(\ref{sturmliouville3})
reduces to the one-dimensional equation~(\ref{sturmliouville1}), as
can be seen by comparing eqs.~(\ref{scaleN}) and~(\ref{shannon1}) after
making the identifications $K=W$ and $R=T$. The fixed-$m$ Slepian
functions $\varphi(\xi)$ can be determined by writing them as the
infinite series 
\be\label{DV}
\varphi(\xi) = m!\,\xi^{m+\frac{1}{2}} \sum_{l=0}^{\infty} 
\fracd{d_{l}\, l!}{(l+m)!} P_{l}^{m0}(1-2\xi^2)
,\qquad 0\le\xi\le 1,
\ee
whereby the $P_{l}^{m0}(x)$ are the Jacobi polynomials~\cite[]{Abramowitz+65}:
\be
P_l^{m0}(x) = \fracd{(l+m)!\,l!}{2^l} \sum_{n=0}^l
\fracd{(x-1)^n(x+1)^{l-n}}{(n+m)!\,(l-n)!^2n!}
.
\ee
By extension to~$\xi>1$ they can also be determined from the
rapidly converging Bessel series  
\be
\varphi(\xi)=\frac{m!}{\gamma}\sum_{l=0}^\infty  
\frac{d_l\,l!}{(l+m)!}
\frac{J_{m+2l+1}(c\hspf \xi)}{\sqrt{c\hspf \xi}},
\quad \xi\in\mathbb{R^+}.
\label{SE}
\ee 
From eq.~(\ref{jlim}) we confirm that~$\varphi(0)=0$. In both cases,
the required fixed-$m$ expansion coefficients~$d_l$ can be determined by 
recursion~\cite[]{Bouwkamp47,Slepian64}, which is, however, rarely
stable. It is instead more practical to determine them as the eigenvectors of
the non-symmetric tridiagonal
matrix~\cite[]{DeVilliers+2003,Shkolnisky2007} that is the spectral
form, with the same eigenvalues, of eq.~(\ref{sturmliouville3}), 
\begin{subequations}\begin{eqnarray}
T_{l+1\,l}&=&-\fracd{c^2\,(m+l+1)^2}{(2l+m+1)(2l+m+2)},\\
T_{ll}&=&\left( 2l + m + \frac{1}{2}\right) \left(
2l+m+\frac{3}{2} \right)\nnr\\
&&{}+\fracd{c^2}{2} \left[1+
\fracd{m^2}{(2l+m)(2l+m+2)} \right],\\
T_{l\,l+1}&=&-\fracd{c^2\,(l+1)^2}{(2l+m+2)(2l+m+3)},
\end{eqnarray}
\end{subequations}
where the parameter~$l$ ranges from~$0$ to
some large value that ensures convergence. Finally, the desired
concentration eigenvalues~$\lambda$ can subsequently be 
obtained by direct integration of eq.~(\ref{BDscaled}), or,
alternatively~\cite[]{Slepian64}, from 
\be\label{eigos}
\lambda=2\gamma^2\sqrt{N\DD},\with\gamma=
\frac{c^{m+1/2}d_0}{2^{m+1}(m+1)!}
\left(\sum_{l=0}^\infty d_l\right)^{-1} 
.
\ee
Neither procedure is particularly accurate for small
values~$\lambda$, though these will be rarely
needed in applications. A more uniformly valid stable recursive
scheme is described elsewhere~\cite[]{Slepian+61,Xiao+2001}.
The overall accuracy of the procedures chosen can be verified with the
aid of the exact formula eq.~(\ref{nsubm3}), as for the fixed-$m$
eigenvalues, 
$
N\DD_m=\sum_{\alpha=1}^{\infty}\lambda_{\alpha}
$.

\begin{figure*}[hbt]\centering
\includegraphics[width=1\textwidth]{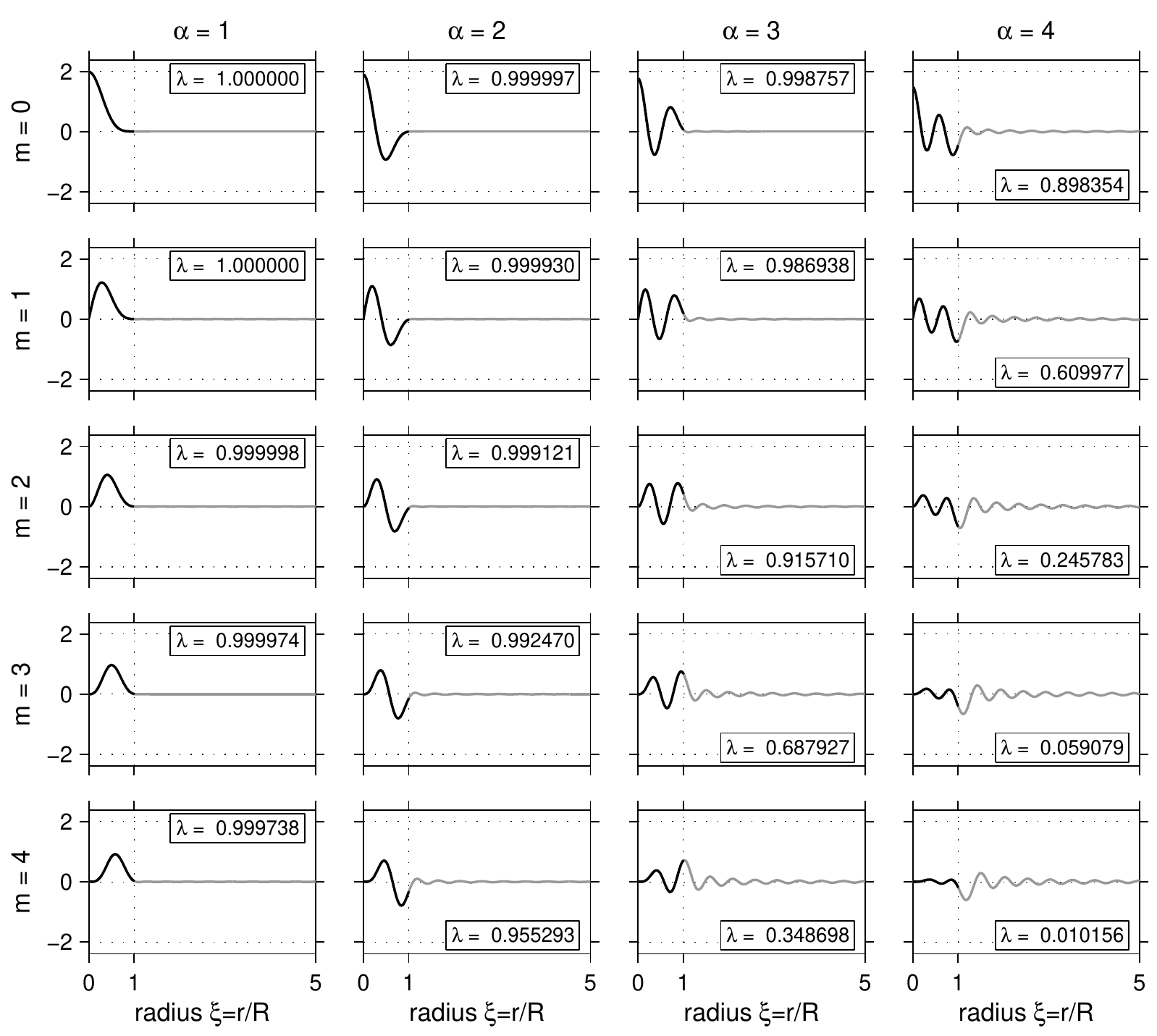}
\caption{\label{swspace2d} Radial dependence of the first four
eigenfunctions of eq.~(\ref{Bfixed}), the bandlimited Slepian
functions $\sfg_\alpha(r)$, $\alpha=1,2,3,4$, of fixed order~$m=0$ (top)
to~$m=4$ (bottom), calculated from eq.~(\ref{SE}) truncated at
$l=84$. The Shannon number~$N\DD=42$ and the bandwidth~$c\approx
13$. Black curves show the concentration within the scaled spatial
interval~$0\le\xi \le 1$ and grey curves show the leakage into the
rest of the real line, which is truncated at~$\xi=5$. Labels show the
eigenvalues, $\lambda_\alpha$, calculated from
eq.~(\ref{eigos}). These express the quality of the spatial
concentration. Spacelimited, spectral-domain Slepian functions are
shown in Fig.~\ref{swspectral2d}.
}
\end{figure*}

\ssec{Numerical examples}\label{numex}

With the eigenfunctions~$\varphi(\xi)$ calculated, by series
expansion~\cite[]{DeVilliers+2003, Shkolnisky2007} via eqs.~(\ref{DV})
or~(\ref{SE}) or even by quadrature~\cite[]{Zhang94} of
eq.~(\ref{GL}), we rescale them following eqs.~(\ref{rescaling})
and~(\ref{Bmapping2}) to the Slepian functions~$\sfg_\alpha(r)$
normalized according to~(\ref{Borthog3}). The sign is arbitrary:
the concentration criterion is a quadratic. 

\begin{figure*}[b]\centering
\includegraphics[width=0.99\textwidth]{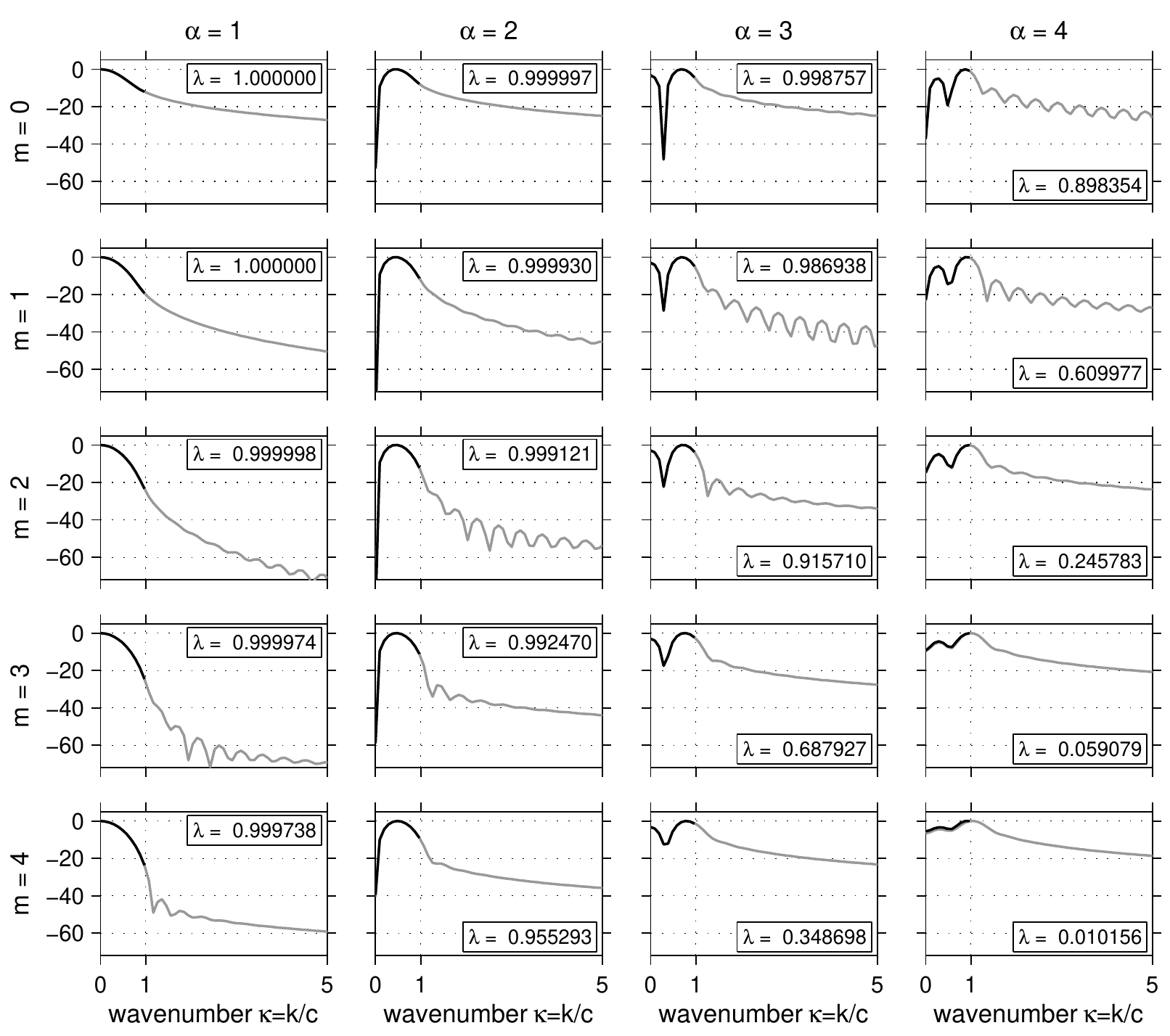}
\caption{\label{swspectral2d}Squared Fourier coefficients~$|H(k)|^2$
of the first four scaled spacelimited eigenfunctions~$\sfh_\alpha(r)$,
$\alpha=1,2,3,4$, of fixed order~$m=0$ (top) to~$m=4$ (bottom),
calculated by discrete Fourier transformation of the functions shown
by the thick black lines in Fig.~\ref{swspace2d}. The Shannon number
is~$N\DD=42$ and the scaled bandwidth~$c\approx 13$.  Black curves show
the power within the scaled wavenumber interval~$0\le \kappa\le 1$ and
grey curves show the power leaked to the rest of the wavenumber axis,
which is truncated at~$\kappa=5$. Values of~$|H(k)|^2$ are in decibel
(dB),  zero at the individual maxima. Labels show the
eigenvalues~$\lambda_\alpha$, which express the quality of the
spectral concentration. The corresponding bandlimited, spatial-domain
eigenfunctions are shown in Fig.~\ref{swspace2d}.}
\end{figure*}

\begin{figure*}[b]\centering 
{\includegraphics[width=0.85\textwidth]{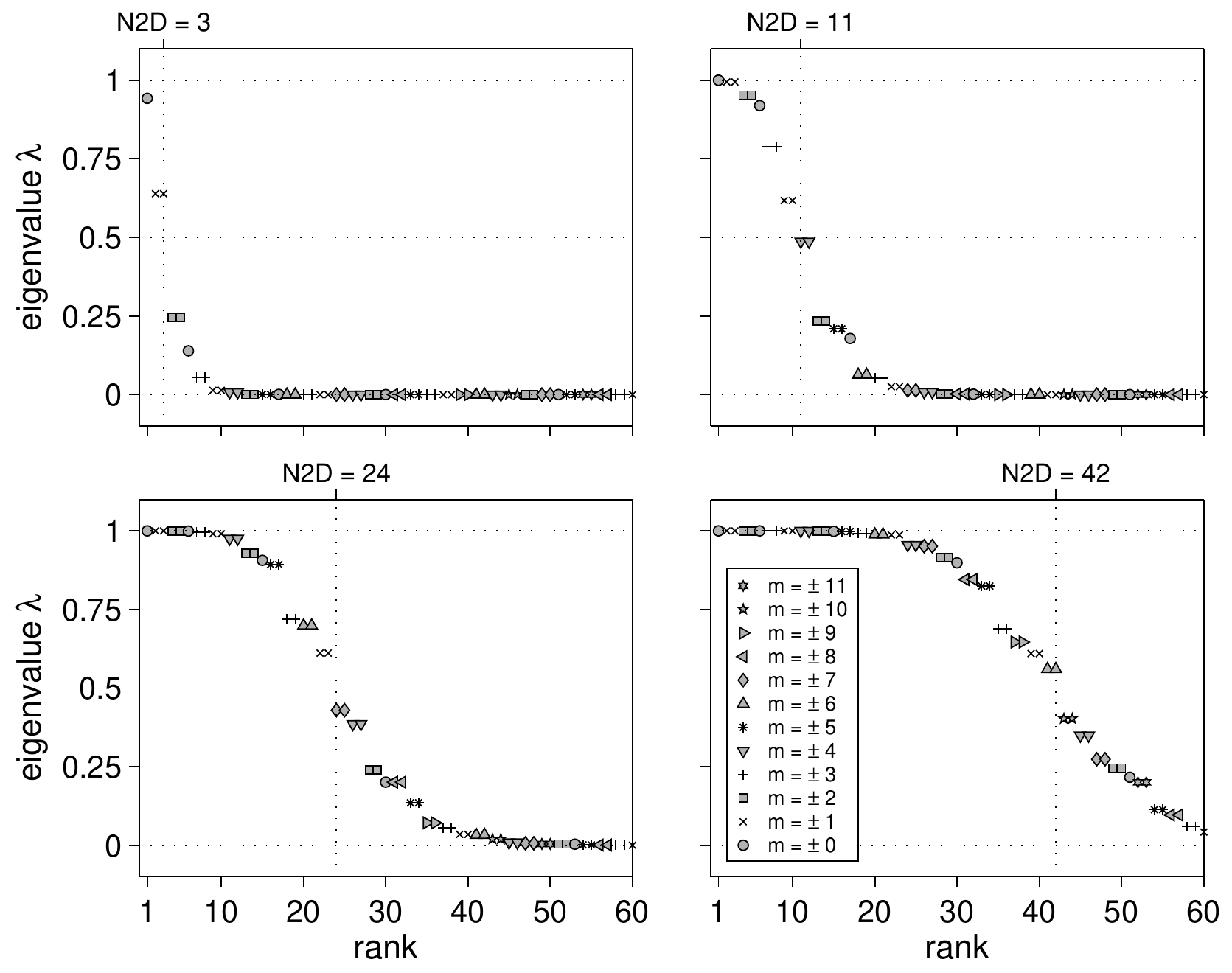}}
\caption{\label{swvals2d}Eigenvalue spectra ($\lambda_{\alpha}$ versus
  rank $\alpha$) of the Cartesian concentration problem, at various
  Shannon numbers~$N\DD$, for disk-shaped regions, calculated via
  eq.~(\ref{eigos}).  The total number of eigenvalues is unbounded;
  only $\lambda_1$ through $\lambda_{60}$ are shown. Different symbols
  are used to plot the various orders $-11\leq m\leq 11$; juxtaposed
  identical symbols are $\pm m$ doublets. Vertical grid lines and top
  labels specify the rounded Shannon numbers $N\DD=3,11,24,42$.}
\end{figure*}

The four most optimally concentrated eigenfunctions
$\sfg_1(r),\sfg_2(r),\sfg_3(r),\sfg_4(r)$, for the orders $0\leq m\leq
4$ are plotted in Fig.~\ref{swspace2d}, scaled for convenience. The
associated eigenvalues $\lambda_1,\lambda_2,\lambda_3,\lambda_4$ are
listed to six-figure accuracy. The planar Shannon number in this
example is $N\DD=42$, hence~$c\approx 13$. It is useful to compare the
behavior of these solutions to those of the spherical concentration
problem~\cite[][Fig.~5.1]{Simons+2006a}, with which they are
asymptotically self-similar~\cite[]{Simons+2007,Simons+2006a}. The first
zeroth-order ``zonal'' ($m=0$) eigenfunction, $\sfg_1(r)$, has no
nodes within the ``cap'' of radius~$R$; the second, $\sfg_2(r)$, has
one node, and so on. The non-zonal ($m>0$) eigenfunctions all vanish
at the origin. The first three~$m=0$, $m=1$ and~$m=2$, and the first
two~$m=3$ and~$m=4$ eigenfunctions are very well concentrated
($\lambda>0.9$). The fourth~$m=3$ and~$m=4$ eigenfunctions exhibit
significant leakage ($\lambda<0.1$).


The squared Fourier coefficients~$|H(k)|^2$ of the four best
concentrated spacelimited eigenfunctions
$\sfh_1(r),\sfh_2(r),\sfh_3(r),\sfh_4(r)$ 
for~$0\leq m\leq 4$ are plotted versus scaled wavenumber~$\kappa=k/c$,
on a decibel scale, in Fig.~\ref{swspectral2d}. Shannon
number, bandwidth, and layout are the same as in
Fig.~\ref{swspace2d}. See elsewhere~\cite[][Fig.~5.2]{Simons+2006a}
to compare with the spherical case.

Once a number of sequences of fixed-order eigenvalues has
been found, they can be resorted to an overall mixed-order
ranking. The total number of significant eigenvalues is then given by
eq.~(\ref{Ntwice}).  In Fig.~\ref{swvals2d} we show the mixed-$m$
eigenvalue spectra for the Shannon numbers 
$N\DD=3,11,24$ and $42$. These show the familiar step
shape~\cite[]{Landau65,Slepian+65}, 
with significant ($\lambda\approx 1$) and insignificant
($\lambda\approx 0$) eigenvalues separated by a narrow transition
band. The Shannon numbers roughly separate the reasonably well
concentrated eigensolutions ($\lambda\ge0.5$) from the more poorly
concentrated ones ($\lambda<0.5$) in all four cases. To compare with
the spherical case, see elsewhere~\cite[][Fig.~5.3]{Simons+2006a}.

\begin{figure*}[tbh]\centering
\includegraphics[width=1\textwidth]{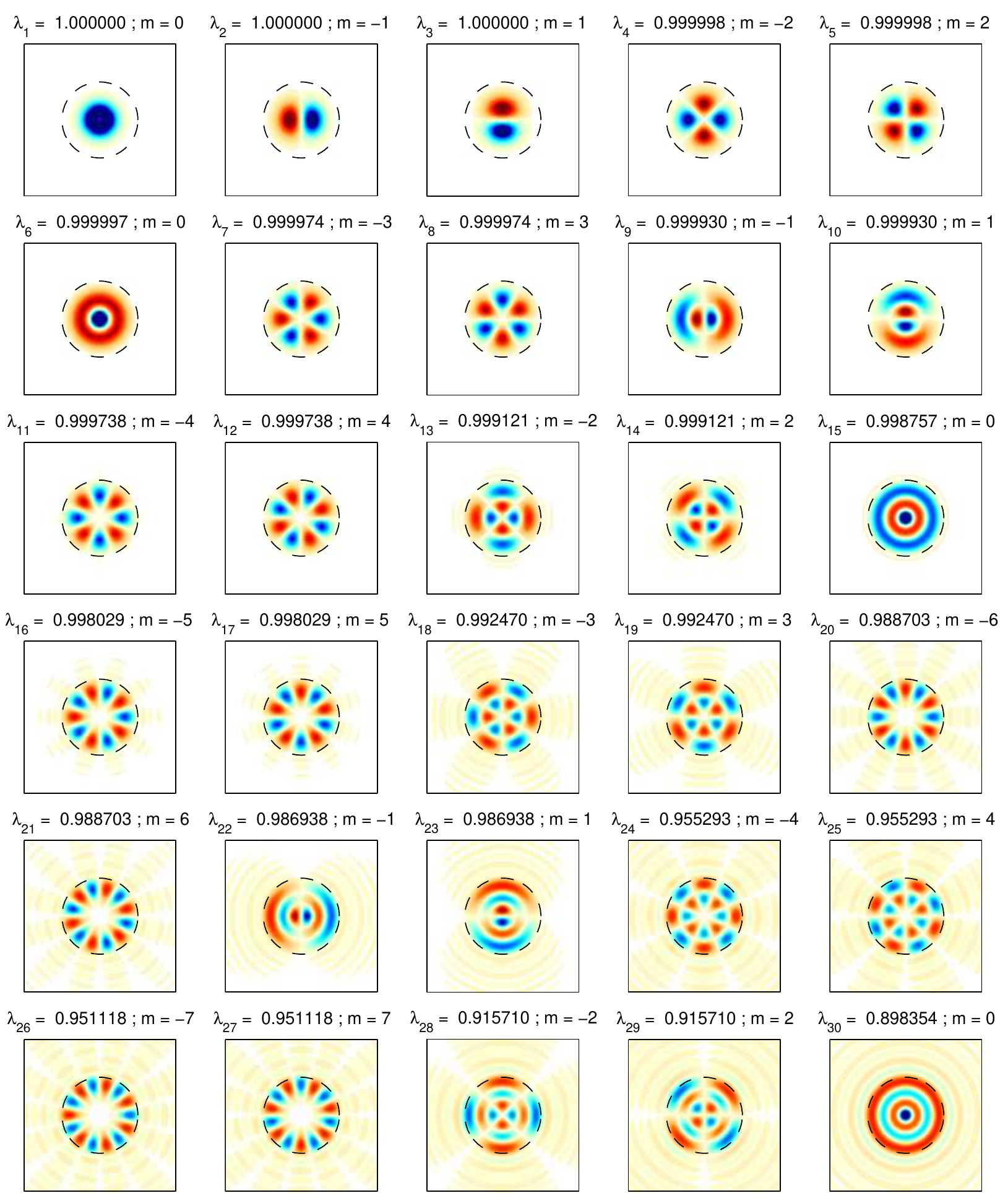}
\caption{\label{swfried2d} Bandlimited eigenfunctions~$g(r,\theta)$
that are optimally concentrated within a disk of radius~$R=1$. Dashed
circle denotes the region boundary. The Shannon number~$N\DD=42$ and
the scaled bandwidth $c\approx 13$. The eigenvalues $\lambda_{\alpha}$
have been sorted into a mixed-order ranking, $\alpha=1\rightarrow 30$, with
the best concentrated eigenfunction plotted on the top left and the
$30$th best on the lower right. Blue (or white in the print version)
is positive and red (black in print) is negative; regions in which the
absolute value is less than one hundredth of the maximum value on the
domain are left white (grey at 50\% in the print version).}  
\end{figure*}
Finally, Fig.~\ref{swfried2d} shows a polar plot of the first 30
eigenfunctions~$g(r,\phi)$ concentrated within a radius $R=1$, defined
by equations~(\ref{Bpolarg}) using the~$\sfg(r)$ shown in
Fig.~\ref{swspace2d}.  The Shannon number is $N^{2D}=42$, as in
Fig.~\ref{swspace2d}. The eigenvalue ranking is mixed-order, as in
Fig.~\ref{swvals2d}, and all degenerate $\sqrt{2}\cos (m\phi),
\sqrt{2}\sin (m\phi)$ doublets are shown. The concentration factors
$1>\lambda\ge 0.8983$ and orders~$m$ of each eigenfunction are
indicated. Blue and red colors (in the online version; in print these
have been converted to grey scales between white and black) represent
positive and negative values, respectively; however, all signs could
be reversed without violating the quadratic concentration
criteria~(\ref{Bnormratio}) and~(\ref{Knormratio}).
Elsewhere~\cite[][Fig.~5.4]{Simons+2006a} these can be compared with
the spherical case. Other calculation methods (see
Sect.~\ref{tocome} below) may yield slight differences between what
should be pairs of eigenvalues for each non-azimuthally symmetric and
thus fixed non-zero order eigenfunction. The mismatch is confined to
the last two quoted digits in a typical
calculation~\cite[compare][Fig.~2]{Simons2010}.

\section{Computational Considerations on General Domains}
\label{tocome}

While the framework of Sect.~\ref{cartesian} is completely
nonrestrictive as to the geometry of the spectral or spatial
concentration regions $\ssK$ or $\ssR$, the
literature treatment has remained focused on the cases with
special symmetry that we discussed in Sect.~\ref{slepsym}, for which
analytic solutions could be found. In what follows we lift all such
restrictions and, numerically, solve the problem of finding the
unknown planar function~$f(\bx)$ that is bandlimited/band-concentrated
to~$\ssK$ and space-concentrated/spacelimited to~$\ssR$,
by satisfying the Fredholm eigenvalue equation
\begin{subequations}
\label{sumup}
\be
\intR \Dxxp
\,f(\bx')\dbx'
=\lambda\hspo f(\bx),
\label{inteqn}\ee\be
\Dxxp=\fnorm\intK e^{i\bk\cdot(\bx-\bx')}\dbk
\label{kernel},
\ee
\end{subequations}
over the domains~$\bx\in\mathbb{R}^2$ (in which case we retrieve the
bandlimited Slepian functions~$g$ of Sect.~\ref{theg}) or $\bx\in 
\ssR$ (when we obtain the spacelimited~$h$ of
Sect.~\ref{theh}).

\subsection{Isotropic spectral response}

Despite the appealing generality of eq.~(\ref{sumup}), applications in
the geophysical or planetary sciences, e.g. in problems of spectral
analysis of data collected on bounded geographical domains, is well
served by Slepian functions whose spectral support is circularly
isotropic. As we saw in Sect.~\ref{circband}, the integral kernel is
then
\be
\Dxxp=\fracd{K\hspo J_1(K\|\bx-\bx'\|)}{2\pi\, \|\bx-\bx'\|}
.
\ee
The standard solution to solve homogeneous Fredholm integral equations
of the second kind is by the Nystr\"om method~\cite[]{Nystrom30,Press+92}. 

\begin{figure*}[b]\centering
\rotatebox{-90}{
\includegraphics[height=\textwidth]{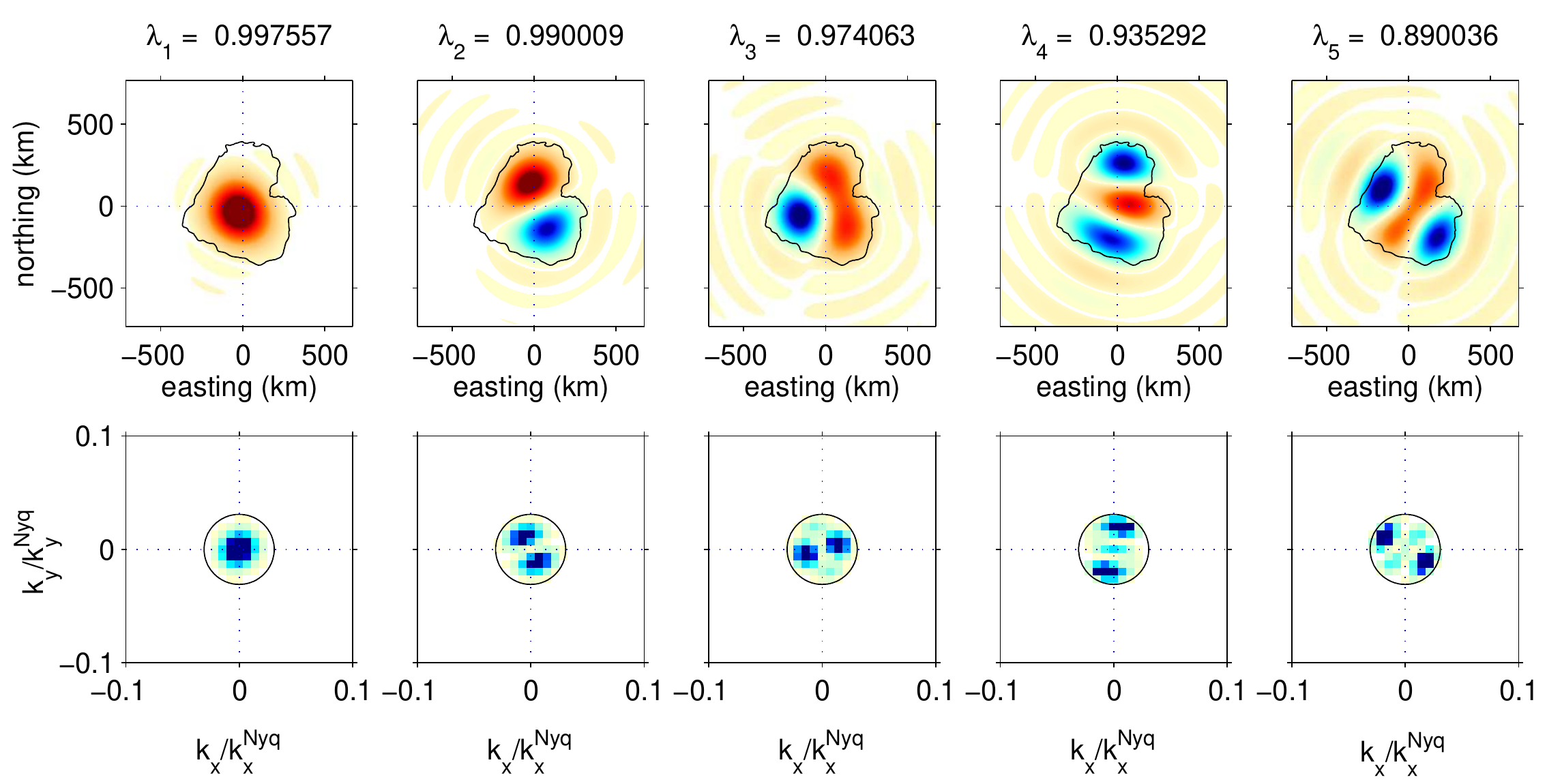}
}
\caption{\label{swregions2d_1}Bandlimited eigenfunctions
$g_1,g_2,\ldots,g_{5}$, calculated by eq.~(\ref{whatwedo}) using
Gauss-Legendre integration with 32 node points in both directions,
optimally concentrated within the Colorado Plateaus, centered on
109.95$^\circ$W 37.01$^\circ$N (near Mexican Hat, Utah) of area
$A\approx334\times 10^3$~km$^2$. The concentration factors
$\lambda_1,\lambda_2,\ldots,\lambda_5$ are indicated; the Shannon
number is $N^{2D}=10$. The top row renders the eigenfunctions in space
on a grid with 5~km resolution in both directions, with the reversible
convention that positive values are blue (white in print) and negative
values red (black in print). The 
spatial concentration region is outlined in black. The bottom row
shows the squared Fourier coefficients $|G(\mbf{k})|^2$ as calculated
from the functions $g(\bx)$ shown, on a wavenumber scale that is
expressed as a fraction of the Nyquist wavenumber. The spectral
limitation region is shown by the black circle at wavenumber
$K=0.0194$~rad/km. All areas for which the absolute value of the
functions plotted is less than one hundredth of the maximum value
attained over the domain are left white (grey at 50\% in the print
version).} 
\end{figure*}

In one dimension, let us write eq.~(\ref{sumup}) in a form
reminiscent of the special case that we encountered previously in
eq.~(\ref{Bfixed}), to which it also applies explicitly, namely    
\be 
\int_{a}^{b} D(x,x') \,f(x') \,dx'= \lambda \hspo f(x)
.
\label{ginteeqn} 
\ee
The left-hand side of this equation is to be discretized by a
quadrature rule; this involves choosing $J$ weights $w_j$ and
abscissas $x_j$ such that, to the desired accuracy,
\be
\int_{a}^{b} D(x,x') \,f(x') \,dx'=
\sum_{j=1}^{J} w_j D(x,x_j)\,f(x_j)
,
\label{discro}
\ee
using which we rewrite eq.~(\ref{ginteeqn}) as
\be
\sum_{j=1}^{J} w_j D(x_i,x_j)\, f(x_j) = \lambda \hspo f(x_i) 
,
\ee
this time evaluating the unknown right hand side~$f(x)$ at the
quadrature points as well. Written in matrix form, we may solve
\be
\mbf{D} \cdot \mbf{W} \cdot \mbf{f} = \lambda \hspo \mbf{f}
,
\label{matrixform}
\ee
identifying the elements of the square kernel $D_{ij}=D(x_i,x_j)$,
the samples of the unknown function as $f_j=f(x_j)$, and the entries of
the diagonal weight matrix $W_{ij}=w_j\delta_{ij}$. 

Restoring the
symmetry of the problem contained in the symmetry
$\mbf{D}=\mbf{D}\Trm$, and provided the weights $w_j$ are positive, we
had rather solve
\be
\mbf{\tilde{W}} \cdot \mbf{D} \cdot \mbf{\tilde{W}} \cdot
\mbf{\tilde{f}} = \lambda \hspo \mbf{\tilde{f}},
\label{syform}
\ee 
using $\tilde{W}_{ij}=\sqrt{w_i}\delta_{ij}$ and~$\mbf{\tilde{f}}=
\mbf{\tilde{W}}\cdot\mbf{f}$ or indeed
$\tilde{f}_i=\sqrt{w_i}f_i$. Having solved this system for the unknown
function~$f(x_i)$ at the integration nodes, we can subsequently
produce either $g(x)$ (anywhere in one-dimensional space) or $h(x)$
(inside the region) using eqs.~(\ref{ginteeqn})--(\ref{discro}), to
the same accuracy.  

\begin{figure*}[b]\centering
{\includegraphics[width=0.75\textwidth]{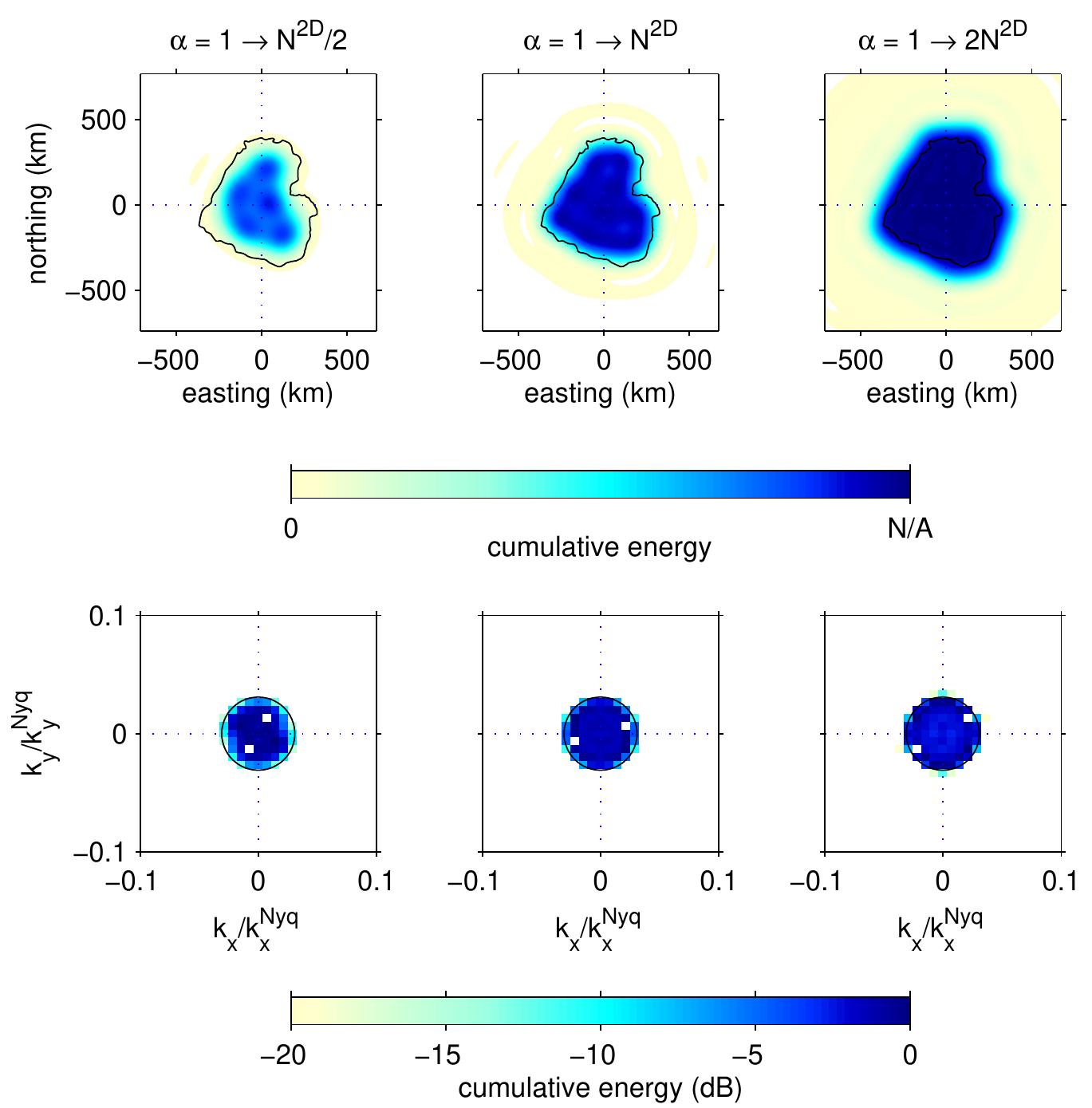}
}
\caption{\label{swconsum2d_1}Cumulative eigenvalue-weighted energy of
the first $N^{2D}/2,N^{2D}$ and $2N^{2D}$ eigenfunctions that are
optimally concentrated within the Colorado Plateaus, with the Shannon
number $N^{2D}=10$. The top row is in the spatial domain with the
concentration region outlined in black. Darkest blue (white in print)
on color bar corresponds to the expected value~(\ref{sumofsq3}) of the
sum, as shown. The bottom row is in the spectral domain with a wavenumber
scale that is expressed as a fraction of the Nyquist wavenumber.  The
spectral limitation region is shown by the black circle. Areas with
absolute value less than one hundredth of the maximum value attained
over the domain are left white (grey at 50\% in the print version).}
\end{figure*}

In two dimensions, the vectors $\bx=(x,y)$ and $\bx'=(x',y')$ denote
points occupying the plane and eq.~(\ref{sumup}) is rewritten as
\be
\intR \Dxxp\,f(\bx')\dbx'=\sum_{k=1}^{K}\sum_{l=1}^{L}
w_k w_l D(\bx;x_k,y_l)\,f(x_k,y_l)
,
\ee
where it is implied that $w_l=w_l(x_k)$ and $y_l=y_l(x_k)$. By
wrapping the indices $k$ and $l$ on the scalar elements $x'$ and $y'$ 
into the vectorial indices $i$ and $j$ used for $\bx$ and $\bx'$ we
arrive at
\be\label{whatwedo}
\sum_{j=1}^{J} w'_j D(\bx_i;\bx_j)\, f(\bx_j) =
 \lambda \hspo f(\bx_i) 
.
\ee
This is identical to eq.~(\ref{syform}) as long as we remember
that the weights $w'_j$ are pairwise products of the
one-dimensional $w_k w_l$ and identifying
$D_{ij}=D(\bx_i;\bx_j)$.

As to the choice of quadrature, a finely meshed Riemann sum might
suffice for some applications~\cite[]{Saito2007,Zhang94}, though we
shall prefer the classical Gauss-Legendre algorithm~\cite[]{Press+92}
for increased accuracy at faster computation speeds. Other options
are available as well~\cite[]{Ramesh+91}. The analytical formulas for the
special cases, e.g. eqs.~(\ref{DV}), (\ref{SE}) and~(\ref{eigos}) can
be used for verification purposes (see the discussion of
Fig.~\ref{swfried2d} in the text), as we have.

Figures~\ref{swregions2d_1} and~\ref{swconsum2d_1} were computed using
eq.~(\ref{whatwedo}) with Gauss-Legendre quadrature, starting from a
splined boundary of the Colorado Plateaus, a physiographic region in
the United States. Fig.~\ref{swregions2d_1} shows the first five
bandlimited Slepian functions, $g_1\rar g_5$, with their eigenvalues,
$\lambda_1\rar \lambda_5$. The space-domain functions~$g(\bx)$, shown
in the top row, are increasingly oscillatory as their concentration
values decrease. At the same time their periodograms~$|G(\mbf{k})|^2$,
shown in the bottom row, remain strictly confined within the chosen
isotropic spectral region of bandwidth~$K$. All together, the Slepian
functions uniformly occupy their respective domains of concentration
or limitation. Fig.~\ref{swconsum2d_1} shows the eigenvalue-weighted
partial sums of the energy of the eigenfunctions, in space (top row)
and spectral space (bottom row). The progressive covering behavior,
consistent with eq.~(\ref{heuristic}), lies at the basis of the
success of the Slepian functions in being used as data tapers for
spectral analysis~\cite[]{Bronez88,Liu+92,Percival+93}: orthogonal over
the domain that they cover, smoothly but rapidly decaying to zero at
the boundary of interest, and with a finite and isotropic spectral
response.
 
\subsection{Arbitrary spectral response}\label{arbit}

We finally treat the situation in which both the spatial and spectral
concentration domains are of arbitrary description --- but with
Hermitian spectral symmetry, see
eqs.~(\ref{Bfourier})--(\ref{BFortho}), if the Slepian functions are to
be real-valued. No analytical results can be expected for such a case,
and we benefit from writing the concentration problem in its most
abstract form. Writing~$\ssQ$ for the operator that acts on a spatial
function~$f(\bx)$ to return its Fourier transform~$F(\bk)$, we
introduce the spatial projection operator
\be\label{oper1}
\ssP f(\bx)=\left\{\barray{ll}
f(\bx) & \mbox{if $\bx\in \ssR$,}\\
0 & \mbox{otherwise},
\earray\right.
\ee
and the spectral projection operator
\be
\ssL F(\bk)=\left\{\barray{ll}
F(\bk) & \mbox{if $\bk\in \ssK$,}\\
0 & \mbox{otherwise}.
\earray\right.
\label{oper2}
\ee
We rewrite the variational equations~(\ref{Bnormratio})
and~(\ref{Knormratio}) in inner-product notation as 
\be\label{inprod}
\lambda=\fracd
{\langle \ssP\ssQi\ssL F, \ssP\ssQi\ssL F\rangle}
{\langle \ssQi\ssL F, \ssQi\ssL F \rangle}=\fracd
{\langle \ssL\ssQ\ssP f, \ssL\ssQ\ssP f\rangle}
{\langle \ssQ\ssP f, \ssQ\ssP f\rangle}
=\mbox{maximum}.
\ee
The associated spectral-domain and spatial-domain eigenvalue
equations are 
\be\label{doit}
\ssL\ssQ\ssP\ssQi\ssL (\ssL F)=\lambda\hspf(\ssL F),\qquad
\ssP\ssQi\ssL\ssQ\ssP (\ssP f)=\lambda\hspf(\ssP f)
,
\ee
where we have made use of the fact that~$\ssP^2=\ssP$
and~$\ssL^2=\ssL$ are self-adjoint, and that~$\ssQ$ and~$\ssQi$ are
each other's adjoints, provided we rewrite eq.~(\ref{Bfourier}) as
a unitary transform. In our notation, the solutions $\ssL F$ yield
the Fourier transforms of the bandlimited functions~$g(\bx)$ of
Sect.~\ref{theg}, while the $\ssP f$ represent the spacelimited
functions~$h(\bx)$ of Sect.~\ref{theh}.

\begin{figure*}[bth]\centering
\rotatebox{-90}{
\includegraphics[height=\textwidth]{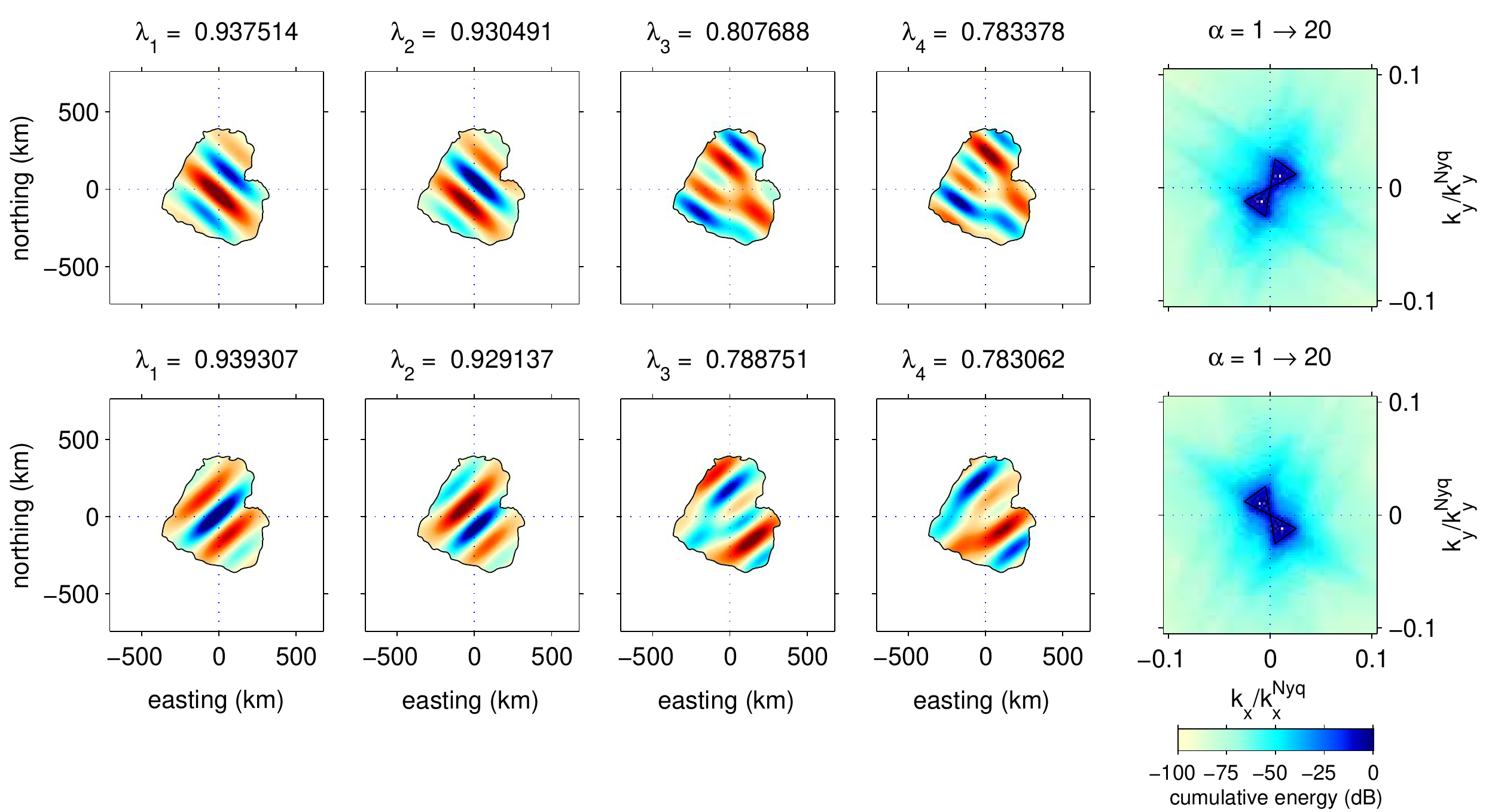}
}
\caption{\label{swsvd_1}Two sets of eigenfunctions
$h_1,h_2,\ldots,h_{4}$ that are strictly spacelimited within the
Colorado Plateaus and optimally concentrated within a wedge-shaped
spectral domain, with their concentration factors
$\lambda_1,\lambda_2,\ldots,\lambda_4$. The last panel on each row
shows the eigenvalue-weighted periodograms of the first~$\alpha=1\rar
20$ eigenfunctions. The functions~$h(\bx)=\ssP f(\bx)$ were calculated
by diagonalization of the spatial-domain operator
equation~(\ref{doit}) on a domain that we took to be three times the
size of the box inscribing the concentration domain outlined in black,
discretized on a 5$\times$5 km$^2$ grid. The spectral concentration
domains are triangles with Hermitian symmetry oriented at~$\pm\pi/6$
from the horizontal in the top and bottom rows, respectively.}
\end{figure*}

For computation in the discrete-discrete case, $\ssP$ and $\ssL$ will
simply be matrices of ones and zeros, and $\ssQ$ and $\ssQi$ will be
the discrete Fourier transform (DFT) matrix and its conjugate
transpose. Both $\ssL\ssQ\ssP\ssQi\ssL$ and $\ssP\ssQi\ssL\ssQ\ssP$
will be Hermitian and typically sparse, hence a number of dedicated
eigenvalues solvers can be used to find the Slepian functions. In our
own \textsc{Matlab} implementation we furthermore take advantage of
defining the operators as ``anonymous'' functions. This greatly
improves the performance of the resulting algorithm, which we tested
by comparison with the methods available for the special cases
discussed in Sections~\ref{realline} and~\ref{slepsym}.  With this we
have complete liberty to design Slepian functions on geographical
domains, including the ability to construct ``steerable'', anisotropic
windows for texture-sensitive analysis that is of interest in theory
\cite[]{Olhede+2009,Vandeville+2008} and practice
\cite[]{Audet+2007,Kirby+2006}. Fig.~\ref{swsvd_1} shows two 
numerical examples of Slepian functions~$h$ spectrally sensitive in
wedge-shaped oriented domains. The top row shows $h_1\rar h_4$ and
their eigenvalues, $\lambda_1\rar \lambda_4$, followed by
eigenvalue-weighted sums of their periodograms,
$\sum_{\alpha=1}^{20}\lambda_\alpha |H_\alpha(\bk)|^2$, on a decibel
scale. These are orientated at $\pi/6$ from the horizontal. The bottom
row shows the equivalent results oriented at $-\pi/6$ from the
horizontal. Both cases were computed by diagonalization
(\texttt{eigs}) of the operator~$\ssP\ssQi\ssL\ssQ\ssP$ whereby~$\ssQ$
is the anonymous function call to the two-dimensional Fast Fourier
Transform (\texttt{fft2}) and~$\ssQi$ its inverse
(\texttt{ifft2}). This immediately returns the~$h(\bx)$ on a
discretized calculation domain in which, once again, the region of the
Colorado Plateaus was embedded. The functions~$h_1\rar h_4$ shown are
the real parts of the eigensolutions associated with the four
eigenvalues~$\lambda_1\rar \lambda_4$ that have the largest real
part. The eigenvalues of these ``generalized prolate spheroidal data
sequences'' (gpss), which need not in principle be defined on
regularly sampled grids~\cite[]{Bronez88}, are sensitive to the
discretization and the size of the computational domains to which the
inner products in eq.~(\ref{inprod}) refer. The slightly different
manner in which these sample the parametrically defined rotated
spectral concentration domains explains the discrepancy of the
eigenvalues between the two cases that are shown in the top and bottom
rows of Fig.~\ref{swsvd_1}.

\section{Conclusion}

Slepian functions are orthogonal families of functions that are all
defined on a common spatial domain, where they are either optimally
concentrated or within which they are exactly limited. At the same
time they are exactly confined within a certain spectral interval, or
maximally concentrated therein~\cite[]{Simons2010}. The joint
optimization of quadratic spatio-spectral concentration is generally
referred to as Slepian's problem, which we encountered as
eq.~(\ref{slepian1}) in one dimension, as eq.~(\ref{normratio}) on the
surface of a sphere, and as eq.~(\ref{Bnormratio}) in two Cartesian
dimensions. Without qualification in one, and under special symmetry
considerations in multiple dimensions, Slepian functions solve
Sturm--Liouville equations, notably eqs.~(\ref{sturmliouville1}),
(\ref{sturmliouville2}), (\ref{sturmliouville22})
and~(\ref{sturmliouville3}). In those cases (though on the sphere only
asymptotically) the solutions depend on the spatial and spectral areas
of concentration only through the Shannon number, which depends on
their product, as defined by eqs.~(\ref{shannon1}), (\ref{shannon2})
and~(\ref{shannon3}). More generally, they solve Fredholm integral
equations in their respective dimensions, as exemplified by
eqs.~(\ref{dude1}), (\ref{firsttimeint}), and~(\ref{Beigen34}).

What unites the Slepian functions is that the kernels of the integral
equations they solve are best thought of in the context of
reproducing-kernel Hilbert
spaces~\cite[]{Amirbekyan+2008b,Nashed+91,Yao67}. Each of
$\Dttp$ in eq.~(\ref{slepian4}), $\Drhrhp$ in  
eq.~(\ref{banddelta}) and $\Dxxp$ in eq.~(\ref{Beigen4}), is a
restriction of the traditional Dirac delta functions, $\delta(t,t')$,
$\delta(\rhat,\rhat')$, and $\delta(\bx,\bx')$, to a region of finite
bandwidth~\cite[]{Simons2010}. While the latter satisfy, for any
square-integrable function~$f$, the relations
\begin{subequations}
\ber 
\label{BFortho0}\label{BFortho00}\label{BFortho000}
f(t)&=& \intinf f(t')\hsp\delta(t,t')\,dt',\\
f(\bx)&=& \intinft f(\bx')\hsp\delta(\bx,\bx')\dbx',\\
f(\rhat)&=& \into f(\rhat')\hsp\delta(\rhat,\rhat')\domg',
\eer
\end{subequations}
the former act on functions~$g$ that belong to appropriately
bandlimited subspaces, as we have defined them in eqs.~(\ref{gband}),
(\ref{bandlg}) and~(\ref{Bgdefn}), in much the same way:
\begin{subequations}
\ber
\label{dttp}\label{dxxp}\label{drrp}
g(t)&=&\intinf g(t')\Dttp\,dt',\\
g(\bx)&=&\intinft g(\bx')\Dxxp\dbx',\\ 
\grh&=&\into \grhp\Drhrhp\domg.
\eer
\end{subequations}
Note that these are \textit{not} equal to eqs.~(\ref{dude1}),
(\ref{firsttimeint}), or~(\ref{Beigen34}). Hence the~$g$ are a basis
for bandlimited processes anywhere on the applicable domain
($\mathbb{R}$, $\mathbb{R}^2$ or the entire spherical
surface~$\Omega$)~\cite[]{Daubechies92,Flandrin98,Freeden+98c,Landau+61,Slepian+61}.
Therein, the Shannon-number best time- or space-concentrated members
allow for sparse, approximate expansions of signals that are
concentrated to the spatial region of choice. Similarly, the infinite
sets of exactly time- or spacelimited (and thus band-concentrated)
functions~$h$ which are the eigenfunctions of eqs.~(\ref{dude1}),
(\ref{firsttimeint}) and~(\ref{Beigen34}) with the domains
appropriately restricted, see eq.~(\ref{Keigen1}) for the
two-dimensional case, are complete bases for square-integrable scalar
functions on the intervals to which they are
confined~\cite[]{Landau+61,Simons+2006a,Slepian+61}.  Expansions of
wideband but spectrally concentrated signals in the small subset of
their Shannon-number most band-concentrated members provide
approximations which are spectrally faithful and constructive as their
numbers grow~\cite[]{Simons2010,Simons+2009b}.

Slepian functions defined on two-dimensional Cartesian domains are of
great utility in the natural sciences. Despite this, no comprehensive
prior treatment was available, beyond that of cases for regions with
advanced symmetry. The present paper has attempted to do right by,
especially, the geosciences, where Fourier-based methods and analyses
are often required but the domains of interest are rarely circles or
rectangles. Much can be learned from those special cases, however,
including from the one-dimensional and spherical situations that we
have reviewed also. Most particularly, they have helped us test and
benchmark the algorithms that we have discussed and are making
available on the Web as part of this publication. In having ended
Sect.~\ref{tocome} with a computational procedure for the design of
planar Cartesian Slepian functions on arbitrarily irregular spatial
and spectral domains, which could be captured in just a few lines, see
eqs.~(\ref{oper1})--(\ref{doit}), and carried out in a handful of
lines of \textsc{Matlab} code, we believe to have returned to a level
of practicality that should appeal to researchers across a wide
spectrum.

\section*{Acknowledgments}
Financial support for this work has been provided in part by the
U.~S.~National Science Foundation under grants EAR-0105387 and
EAR-0710860 to FJS. DVW was also supported by a first-year Graduate
Fellowship from Princeton University. Our Princeton colleagues Tony
Dahlen (1942--2007), Jeremy Goodman (Astrophysical Sciences) and
Eugene Brevdo (Electrical Engineering) each provided a key insight
that helped shape the manuscript. We thank Laura Larsen-Strecker for
help with the identification of geologic provinces, and Ignace Loris,
Volker Michel and Mark Wieczorek for discussions. The comments by the
editor, Willi Freeden, and two anonymous reviewers were all much
appreciated and helped improve the manuscript. All computer code is
posted on \url{www.frederik.net}. 

\bibliography{/home/fjsimons/BIBLIO/bib}

\begin{thebibliography}{130}
\providecommand{\natexlab}[1]{#1}
\expandafter\ifx\csname urlstyle\endcsname\relax
  \providecommand{\doi}[1]{doi:\discretionary{}{}{}#1}\else
  \providecommand{\doi}{doi:\discretionary{}{}{}\begingroup
  \urlstyle{rm}\Url}\fi

\bibitem[{\textit{Abramowitz and Stegun}(1965)}]{Abramowitz+65}
Abramowitz, M., and I.~A. Stegun, \textit{Handbook of Mathematical Functions},
  Dover, New York, 1965.

\bibitem[{\textit{Albertella and Sacerdote}(2001)}]{Albertella+2001}
Albertella, A., and F.~Sacerdote, Using {S}lepian functions for local geodetic
  computations, \textit{Boll.~Geod.~Sc.~Aff.}, \textit{60}(1), 1--14, 2001.

\bibitem[{\textit{Albertella et~al.}(1999)\textit{Albertella, Sans{\`o}, and
  Sneeuw}}]{Albertella+99}
Albertella, A., F.~Sans{\`o}, and N.~Sneeuw, Band-limited functions on a
  bounded spherical domain: the {S}lepian problem on the sphere,
  \textit{J.~Geodesy}, \textit{73}, 436--447, 1999.

\bibitem[{\textit{Amirbekyan et~al.}(2008)\textit{Amirbekyan, Michel, and
  Simons}}]{Amirbekyan+2008b}
Amirbekyan, A., V.~Michel, and F.~J. Simons, Parameterizing surface-wave
  tomopgraphic models with harmonic spherical splines,
  \textit{Geophys.~J.~Int.}, \textit{174}(2), 617, doi:
  10.1111/j.1365--246X.2008.03,809.x, 2008.

\bibitem[{\textit{Aronszajn}(1950)}]{Aronszajn50}
Aronszajn, N., Theory of reproducing kernels, \textit{Trans.~Am.~Math.~Soc.},
  \textit{68}(3), 337--404, 1950.

\bibitem[{\textit{Audet and Mareschal}(2007)}]{Audet+2007}
Audet, P., and J.-C. Mareschal, Wavelet analysis of the coherence between
  {B}ouguer gravity and topography: application to the elastic thickness
  anisotropy in the {C}anadian {S}hield, \textit{Geophys.~J.~Int.},
  \textit{168}, 287--298, doi: 10.1111/j.1365--246X.2006.03,231.x, 2007.

\bibitem[{\textit{Bell et~al.}(1993)\textit{Bell, Percival, and
  Walden}}]{Bell+93}
Bell, B., D.~B. Percival, and A.~T. Walden, Calculating {T}homson's spectral
  multitapers by inverse iteration, \textit{J.~Comp.~Graph.~Stat.},
  \textit{2}(1), 119--130, 1993.

\bibitem[{\textit{Bertero et~al.}(1985{\natexlab{a}})\textit{Bertero, de~Mol,
  and Pike}}]{Bertero+85}
Bertero, M., C.~de~Mol, and E.~R. Pike, Linear inverse problems with discrete
  data.~{I}.~{G}eneral formulation and singular system analysis,
  \textit{Inv.~Probl.}, \textit{1}, 301--330, doi: 10.1088/0266--5611/1/4/004,
  1985{\natexlab{a}}.

\bibitem[{\textit{Bertero et~al.}(1985{\natexlab{b}})\textit{Bertero, de~Mol,
  and Pike}}]{Bertero+88}
Bertero, M., C.~de~Mol, and E.~R. Pike, Linear inverse problems with discrete
  data.~{II}.~{S}tability and regularisation, \textit{Inv.~Probl.}, \textit{1},
  301--330, doi: 10.1088/0266--5611/1/4/004, 1985{\natexlab{b}}.

\bibitem[{\textit{Beylkin and Monz\'on}(2002)}]{Beylkin+2002}
Beylkin, G., and L.~Monz\'on, On generalized {G}aussian quadratures for
  exponentials and their applications, \textit{Appl.~Comput.~Harmon.~Anal.},
  \textit{12}, 332--372, doi: 10.1006/acha.2002.0380, 2002.

\bibitem[{\textit{Beylkin and Sandberg}(2005)}]{Beylkin+2005}
Beylkin, G., and K.~Sandberg, Wave propagation using bases for bandlimited
  functions, \textit{Wave Motion}, \textit{41}(3), 263--291, 2005.

\bibitem[{\textit{Borcea et~al.}(2008)\textit{Borcea, Papanicolaou, and
  Vasquez}}]{Borcea+2008}
Borcea, L., G.~Papanicolaou, and F.~G. Vasquez, Edge illumination and imaging
  of extended reflectors, \textit{SIAM J.~Imag.~Sci.}, \textit{1}(1), 75--114,
  doi:10.1137/07069,290X, 2008.

\bibitem[{\textit{Bouwkamp}(1947)}]{Bouwkamp47}
Bouwkamp, C.~J., On spheroidal wave functions of order zero,
  \textit{J.~Math.~Phys.}, \textit{26}, 79--92, 1947.

\bibitem[{\textit{Boyd}(2003)}]{Boyd2003}
Boyd, J.~P., Approximation of an analytic function on a finite real interval by
  a bandlimited function and conjectures on properties of prolate spheroidal
  functions, \textit{Appl.~Comput.~Harmon.~Anal.}, \textit{15}(2), 168--176,
  2003.

\bibitem[{\textit{Boyd}(2004)}]{Boyd2004}
Boyd, J.~P., Prolate spheroidal wavefunctions as an alternative to {C}hebyshev
  and {L}egendre polynomials for spectral element and pseudospectral
  algorithms, \textit{J.~Comput.~Phys.}, \textit{199}(2), 688--716, 2004.

\bibitem[{\textit{Brander and DeFacio}(1986)}]{Brander+86}
Brander, O., and B.~DeFacio, A generalisation of {S}lepian's solution for the
  singular value decomposition of filtered {F}ourier transforms,
  \textit{Inv.~Probl.}, \textit{2}, L9--L14, 1986.

\bibitem[{\textit{Bronez}(1988)}]{Bronez88}
Bronez, T.~P., Spectral estimation of irregularly sampled multidimensional
  processes by generalized prolate spheroidal sequences, \textit{IEEE
  Trans.~Acoust.~Speech~Signal~Process.}, \textit{36}(12), 1862--1873, 1988.

\bibitem[{\textit{Chambodut et~al.}(2005)\textit{Chambodut, Panet, Mandea,
  Diament, Holschneider, and Jamet}}]{Chambodut+2005}
Chambodut, A., I.~Panet, M.~Mandea, M.~Diament, M.~Holschneider, and O.~Jamet,
  Wavelet frames: an alternative to spherical harmonic representation of
  potential fields, \textit{Geophys.~J.~Int.}, \textit{163}(3), 875--899, 2005.

\bibitem[{\textit{Chen et~al.}(2005)\textit{Chen, Gottlieb, and
  Hesthaven}}]{Chen+2005}
Chen, Q.~Y., D.~Gottlieb, and J.~S. Hesthaven, Spectral methods based on
  prolate spheroidal wave functions for hyperbolic {PDE}s, \textit{Wave
  Motion}, \textit{43}(5), 1912--1933, 2005.

\bibitem[{\textit{Coifman and Lafon}(2006)}]{Coifman+2006}
Coifman, R.~R., and S.~Lafon, Geometric harmonics: {A} novel tool for
  multiscale out-of-sample extension of empirical functions,
  \textit{Appl.~Comput.~Harmon.~Anal.}, \textit{21}, 31--52, doi:
  10.1016/j.acha.2005.07.005, 2006.

\bibitem[{\textit{Dahlen and Simons}(2008)}]{Dahlen+2008}
Dahlen, F.~A., and F.~J. Simons, Spectral estimation on a sphere in geophysics
  and cosmology, \textit{Geophys.~J.~Int.}, \textit{174}, 774--807, doi:
  10.1111/j.1365--246X.2008.03,854.x, 2008.

\bibitem[{\textit{Dahlen and Tromp}(1998)}]{Dahlen+98}
Dahlen, F.~A., and J.~Tromp, \textit{Theoretical Global Seismology}, Princeton
  Univ.~Press, Princeton, N.~J., 1998.

\bibitem[{\textit{Daubechies}(1988)}]{Daubechies88a}
Daubechies, I., Time-frequency localization operators: {A} geometric phase
  space approach, \textit{IEEE Trans.~Inform.~Theory}, \textit{34}, 605--612,
  1988.

\bibitem[{\textit{Daubechies}(1992)}]{Daubechies92}
Daubechies, I., \textit{Ten Lectures on Wavelets}, \textit{CBMS-NSF Regional
  Conference Series in Applied Mathematics}, vol.~61, Society for Industrial \&
  Applied Mathematics, Philadelphia, Penn., 1992.

\bibitem[{\textit{Daubechies and Paul}(1988)}]{Daubechies+88}
Daubechies, I., and T.~Paul, Time-frequency localisation operators --- {A}
  geometric phase space approach: {II}. {T}he use of dilations,
  \textit{Inv.~Probl.}, \textit{4}(3), 661--680, 1988.

\bibitem[{\textit{de~Villiers et~al.}(2001)\textit{de~Villiers, Marchaud, and
  Pike}}]{DeVilliers+2001}
de~Villiers, G.~D., F.~B.~T. Marchaud, and E.~R. Pike, Generalized {G}aussian
  quadrature applied to an inverse problem in antenna theory, \textit{Inverse
  Problems}, \textit{17}, 1163--1179, 2001.

\bibitem[{\textit{de~Villiers et~al.}(2003)\textit{de~Villiers, Marchaud, and
  Pike}}]{DeVilliers+2003}
de~Villiers, G.~D., F.~B.~T. Marchaud, and E.~R. Pike, Generalized {G}aussian
  quadrature applied to an inverse problem in antenna theory: {II.} {T}he
  two-dimensional case with circular symmetry, \textit{Inverse Problems},
  \textit{19}, 755--778, 2003.

\bibitem[{\textit{Delsarte et~al.}(1985)\textit{Delsarte, Janssen, and
  Vries}}]{Delsarte+85}
Delsarte, P., A.~J. E.~M. Janssen, and L.~B. Vries, Discrete prolate spheroidal
  wave functions and interpolation, \textit{SIAM J.~Appl.~Math.},
  \textit{45}(4), 641--650, 1985.

\bibitem[{\textit{Donoho and Stark}(1989)}]{Donoho+89}
Donoho, D.~L., and P.~B. Stark, Uncertainty principles and signal recovery,
  \textit{SIAM J.~Appl.~Math.}, \textit{49}(3), 906--931, 1989.

\bibitem[{\textit{Edmonds}(1996)}]{Edmonds96}
Edmonds, A.~R., \textit{Angular Momentum in Quantum Mechanics}, Princeton
  Univ.~Press, Princeton, N.J., 1996.

\bibitem[{\textit{Evans et~al.}(2010)\textit{Evans, Andrews-Hanna, and
  Zuber}}]{Evans+2010}
Evans, A.~J., J.~C. Andrews-Hanna, and M.~T. Zuber, Geophysical limitations on
  the erosion history within {A}rabia {T}erra, \textit{J.~Geophys.~Res.},
  \textit{115}, E05,007, doi: 10.1029/2009JE003,469, 2010.

\bibitem[{\textit{Fa\"{y} et~al.}(2008)\textit{Fa\"{y}, Guilloux, Betoule,
  Cardoso, Delabrouille, and Jeune}}]{Fay+2008}
Fa\"{y}, G., F.~Guilloux, M.~Betoule, J.-F. Cardoso, J.~Delabrouille, and M.~L.
  Jeune, {CMB} power spectrum estimation using wavelets, \textit{Phys.~Rev.\
  D}, \textit{78}, 083,013, doi: 10.1103/PhysRevD.78.083,013, 2008.

\bibitem[{\textit{Fengler et~al.}(2007)\textit{Fengler, Freeden, Kohlhaas,
  Michel, and Peters}}]{Fengler+2006a}
Fengler, M.~J., W.~Freeden, A.~Kohlhaas, V.~Michel, and T.~Peters, Wavelet
  modeling of regional and temporal variations of the earth's gravitational
  potential observed by {GRACE}, \textit{J.~Geodesy}, \textit{81}(1), 5--15,
  doi: 10.1007/s00,190--006--0040--1, 2007.

\bibitem[{\textit{Flandrin}(1998)}]{Flandrin98}
Flandrin, P., \textit{Temps-{F}r{\'e}quence}, 2 ed., Herm{\`e}s, Paris, 1998.

\bibitem[{\textit{Freeden and Michel}(1999)}]{Freeden+99}
Freeden, W., and V.~Michel, Constructive approximation and numerical methods in
  geodetic research today -- an attempt at a categorization based on an
  uncertainty principle, \textit{J.~Geodesy}, \textit{73}(9), 452--465, 1999.

\bibitem[{\textit{Freeden and Windheuser}(1997)}]{Freeden+97}
Freeden, W., and U.~Windheuser, Combined spherical harmonic and wavelet
  expansion --- {A} future concept in {E}arth's gravitational determination,
  \textit{Appl.~Comput.~Harmon.~Anal.}, \textit{4}, 1--37, 1997.

\bibitem[{\textit{Freeden et~al.}(1998)\textit{Freeden, Gervens, and
  Schreiner}}]{Freeden+98c}
Freeden, W., T.~Gervens, and M.~Schreiner, \textit{Constructive Approximation
  {o}n {t}he Sphere}, Clarendon Press, Oxford, UK, 1998.

\bibitem[{\textit{Golub and van Loan}(1989)}]{Golub+89}
Golub, G.~H., and C.~F. van Loan, \textit{Matrix Computations}, 2nd ed., Johns
  Hopkins Univ.~Press, Baltimore, Md., 1989.

\bibitem[{\textit{Gradshteyn and Ryzhik}(2000)}]{Gradshteyn+2000}
Gradshteyn, I.~S., and I.~M. Ryzhik, \textit{Tables of Integrals, Series, and
  Products}, 6 ed., Academic Press, San Diego, Calif., 2000.

\bibitem[{\textit{Gr{\"u}nbaum}(1981)}]{Grunbaum81a}
Gr{\"u}nbaum, F.~A., Eigenvectors of a {T}oeplitz matrix: discrete version of
  the prolate spheroidal wave functions, \textit{SIAM J.~Alg.~Disc.~Meth.},
  \textit{2}(2), 136--141, 1981.

\bibitem[{\textit{Gr{\"u}nbaum et~al.}(1982)\textit{Gr{\"u}nbaum, Longhi, and
  Perlstadt}}]{Grunbaum+82}
Gr{\"u}nbaum, F.~A., L.~Longhi, and M.~Perlstadt, Differential operators
  commuting with finite convolution integral operators: some non-{A}belian
  examples, \textit{SIAM J.~Appl.~Math.}, \textit{42}(5), 941--955, 1982.

\bibitem[{\textit{Hall and Mitchell}(2002)}]{Hall+2002}
Hall, B.~C., and J.~J. Mitchell, Coherent states on spheres,
  \textit{J.~Math.~Phys.}, \textit{43}(3), 1211--1236, 2002.

\bibitem[{\textit{Han}(2008)}]{Han2008}
Han, S.-C., Improved regional gravity fields on the {M}oon from {L}unar
  {P}rospector tracking data by means of localized spherical harmonic
  functions, \textit{J.~Geophys.~Res.}, \textit{113}, E11,012,
  doi:10.1029/2008JE003,166, 2008.

\bibitem[{\textit{Han and Ditmar}(2007)}]{Han+2007}
Han, S.-C., and P.~Ditmar, Localized spectral analysis of global satellite
  gravity fields for recovering time-variable mass redistributions,
  \textit{J.~Geodesy}, \textit{82}(7), 423--430, doi:
  10.1007/s00,190--007--0194--5, 2007.

\bibitem[{\textit{Han and Simons}(2008)}]{Han+2008a}
Han, S.-C., and F.~J. Simons, Spatiospectral localization of global
  geopotential fields from the {G}ravity {R}ecovery and {C}limate {E}xperiment
  {GRACE} reveals the coseismic gravity change owing to the 2004
  {S}umatra-{A}ndaman earthquake, \textit{J.~Geophys.~Res.}, \textit{113},
  B01,405, doi: 10.1029/2007JB004,927, 2008.

\bibitem[{\textit{Han et~al.}(2008{\natexlab{a}})\textit{Han, Rowlands,
  Luthcke, and Lemoine}}]{Han+2008b}
Han, S.-C., D.~D. Rowlands, S.~B. Luthcke, and F.~G. Lemoine, Localized
  analysis of satellite tracking data for studying time-variable {E}arth's
  gravity fields, \textit{J.~Geophys.~Res.}, \textit{113}, B06,401, doi:
  10.1029/2007JB005,218, 2008{\natexlab{a}}.

\bibitem[{\textit{Han et~al.}(2008{\natexlab{b}})\textit{Han, Sauber, Luthcke,
  Ji, and Pollitz}}]{Han+2008c}
Han, S.-C., J.~Sauber, S.~B. Luthcke, C.~Ji, and F.~F. Pollitz, Implications of
  postseismic gravity change following the great 2004 {S}umatra-{A}ndaman
  earthquake from the regional harmonic analysis of {GRACE} inter-satellite
  tracking data, \textit{J.~Geophys.~Res.}, \textit{113}, B11,413, doi:
  10.1029/2008JB005,705, 2008{\natexlab{b}}.

\bibitem[{\textit{Han et~al.}(2009)\textit{Han, Mazarico, and
  Lemoine}}]{Han+2009}
Han, S.-C., E.~Mazarico, and F.~G. Lemoine, Improved nearside gravity field of
  the {M}oon by localizing the power law constraint,
  \textit{Geophys.~Res.~Lett.}, \textit{36}, L11,203,
  doi:10.1029/2009GL038,556, 2009.

\bibitem[{\textit{Hanssen}(1997)}]{Hanssen97}
Hanssen, A., Multidimensional multitaper spectral estimation, \textit{Signal
  Process.}, \textit{58}, 327--332, 1997.

\bibitem[{\textit{Harig et~al.}(2010)\textit{Harig, Zhong, and
  Simons}}]{Harig+2010}
Harig, C., S.~Zhong, and F.~J. Simons, Constraints on upper-mantle viscosity
  inferred from the flow-induced pressure gradient across a continental keel,
  \textit{Geochem.~Geophys.~Geosys.}, \textit{11}(6), Q06,004, doi:
  10.1029/2010GC003,038, 2010.

\bibitem[{\textit{Holschneider et~al.}(2003)\textit{Holschneider, Chambodut,
  and Mandea}}]{Holschneider+2003}
Holschneider, M., A.~Chambodut, and M.~Mandea, From global to regional analysis
  of the magnetic field on the sphere using wavelet frames, \textit{Phys.~Earth
  Planet.~Inter.}, \textit{135}, 107--124, 2003.

\bibitem[{\textit{Jackson et~al.}(1991)\textit{Jackson, Meyer, Nishimura, and
  Macovski}}]{Jackson+91}
Jackson, J.~I., C.~H. Meyer, D.~G. Nishimura, and A.~Macovski, Selection of a
  convolution function for {F}ourier inversion using gridding, \textit{IEEE
  Trans.~Med.~Imag.}, \textit{10}(3), 473--478, 1991.

\bibitem[{\textit{Jeffreys and Jeffreys}(1988)}]{Jeffreys+88}
Jeffreys, H., and B.~S. Jeffreys, \textit{Methods of Mathematical Physics}, 3
  ed., Cambridge Univ.~Press, Cambridge, UK, 1988.

\bibitem[{\textit{Karoui and Moumni}(2008)}]{Karoui+2008}
Karoui, A., and T.~Moumni, New efficient methods of computing the prolate
  spheroidal wave functions and their corresponding eigenvalues,
  \textit{Appl.~Comput.~Harmon.~Anal.}, \textit{24}(3), 269--289, 2008.

\bibitem[{\textit{Kennedy et~al.}(2008)\textit{Kennedy, Zhang, and
  Abhayapala}}]{Kennedy+2008}
Kennedy, R.~A., W.~Zhang, and T.~D. Abhayapala, Spherical harmonic analysis and
  model-limited extrapolation on the sphere: {I}ntegral equation formulation,
  in \textit{Proc.~IEEE Int.~Conf.~Signal Process.~Comm.~Syst.}, pp. 1--6, doi:
  10.1109/ICSPCS.2008.4813,702, IEEE, 2008.

\bibitem[{\textit{Khare and George}(2003)}]{Khare+2003}
Khare, K., and N.~George, Sampling theory approach to prolate spheroidal
  wavefunctions, \textit{J.~Phys.~A:\ Math.~Gen.}, \textit{36}, 10,011--10,021,
  2003.

\bibitem[{\textit{Kido et~al.}(2003)\textit{Kido, Yuen, and
  Vincent}}]{Kido+2003}
Kido, M., D.~A. Yuen, and A.~P. Vincent, Continuous wavelet-like filter for a
  spherical surface and its application to localized admittance function on
  {M}ars, \textit{Phys.~Earth Planet.~Inter.}, \textit{135}, 1--14, 2003.

\bibitem[{\textit{Kirby and Swain}(2006)}]{Kirby+2006}
Kirby, J.~F., and C.~J. Swain, Mapping the mechanical anisotropy of the
  lithosphere using a {2D} wavelet coherence, and its application to
  {A}ustralia, \textit{Phys.~Earth Planet.~Inter.}, \textit{158}(2--4),
  122--138, doi: 10.1016/j.pepi.2006.03.022, 2006.

\bibitem[{\textit{Kowalski and Rembieli{\'n}ski}(2000)}]{Kowalski+2000}
Kowalski, K., and J.~Rembieli{\'n}ski, Quantum mechanics on a sphere and
  coherent states, \textit{J.~Phys.~A:\ Math.~Gen.}, \textit{33}, 6035--6048,
  2000.

\bibitem[{\textit{Lai et~al.}(2009)\textit{Lai, Shum, Baramidze, and
  Wenston}}]{Lai+2009}
Lai, M.~J., C.~K. Shum, V.~Baramidze, and P.~Wenston, Triangulated spherical
  splines for geopotential reconstruction, \textit{J.~Geod.}, \textit{83},
  695--708, doi: 10.1007/s00,190--008--0283--0, 2009.

\bibitem[{\textit{Landau}(1965)}]{Landau65}
Landau, H.~J., On the eigenvalue behavior of certain convolution equations,
  \textit{Trans.~Am.~Math.~Soc.}, \textit{115}, 242--256, 1965.

\bibitem[{\textit{Landau and Pollak}(1961)}]{Landau+61}
Landau, H.~J., and H.~O. Pollak, Prolate spheroidal wave functions, {F}ourier
  analysis and uncertainty --- {II}, \textit{Bell Syst.~Tech.~J.},
  \textit{40}(1), 65--84, 1961.

\bibitem[{\textit{Landau and Pollak}(1962)}]{Landau+62}
Landau, H.~J., and H.~O. Pollak, Prolate spheroidal wave functions, {F}ourier
  analysis and uncertainty --- {III}: {T}he dimension of the space of
  essentially time- and band-limited signals, \textit{Bell Syst.~Tech.~J.},
  \textit{41}(4), 1295--1336, 1962.

\bibitem[{\textit{Lilly and Park}(1995)}]{Lilly+95}
Lilly, J.~M., and J.~Park, Multiwavelet spectral and polarization analyses of
  seismic records, \textit{Geophys.~J.~Int.}, \textit{122}, 1001--1021, 1995.

\bibitem[{\textit{Lindquist et~al.}(2006)\textit{Lindquist, Zhang, Glover,
  Shepp, and Yang}}]{Lindquist+2006}
Lindquist, M.~A., C.~H. Zhang, G.~Glover, L.~Shepp, and Q.~X. Yang, A
  generalization of the two-dimensional prolate spheroidal wave function method
  for nonrectilinear {MRI} data acquisition methods, \textit{IEEE
  Trans.~Image~Proc.}, \textit{15}(9), 2792--2804, doi:
  10.1109/TIP.2006.877,314, 2006.

\bibitem[{\textit{Liu and van Veen}(1992)}]{Liu+92}
Liu, T.-C., and B.~D. van Veen, Multiple window based minimum variance spectrum
  estimation for multidimensional random fields, \textit{IEEE
  Trans.~Signal~Process.}, \textit{40}(3), 578--589, doi: 10.1109/78.120,801,
  1992.

\bibitem[{\textit{Ma et~al.}(1996)\textit{Ma, Rokhlin, and Wandzura}}]{Ma+96}
Ma, J., V.~Rokhlin, and S.~Wandzura, Generalized {G}aussian quadrature rules
  for systems of arbitrary functions, \textit{SIAM J.~Numer.~Anal.},
  \textit{33}(3), 971--996, 1996.

\bibitem[{\textit{Mallat}(1998)}]{Mallat98}
Mallat, S., \textit{A Wavelet Tour {o}f Signal Processing}, Academic Press, San
  Diego, Calif., 1998.

\bibitem[{\textit{Maniar and Mitra}(2005)}]{Maniar+2005}
Maniar, H., and P.~P. Mitra, The concentration problem for vector fields,
  \textit{Int.~J.~Bioelectromagn.}, \textit{7}(1), 142--145, 2005.

\bibitem[{\textit{Marinucci et~al.}(2008)\textit{Marinucci, Pietrobon, Balbi,
  Baldi, Cabella, Kerkyacharian, Natoli, Picard, and
  Vittorio}}]{Marinucci+2008}
Marinucci, D., D.~Pietrobon, A.~Balbi, P.~Baldi, P.~Cabella, G.~Kerkyacharian,
  P.~Natoli, Picard, and N.~Vittorio, Spherical needlets for cosmic microwave
  background data analysis, \textit{Mon.~Not.~R.~Astron.~Soc}, \textit{383}(2),
  539--545, doi: 10.1111/j.1365--2966.2007.12,550.x, 2008.

\bibitem[{\textit{McEwen et~al.}(2007)\textit{McEwen, Hobson, Mortlock, and
  Lasenby}}]{McEwen+2007}
McEwen, J.~D., M.~P. Hobson, D.~J. Mortlock, and A.~N. Lasenby, Fast
  directional continuous spherical wavelet transform algorithms, \textit{IEEE
  Trans.~Signal~Process.}, \textit{55}(2), 520--529, 2007.

\bibitem[{\textit{Michel and Wolf}(2008)}]{Michel+2008}
Michel, V., and K.~Wolf, Numerical aspects of a spline-based multiresolution
  recovery of the harmonic mass density out of gravity functionals,
  \textit{Geophys.~J.~Int.}, \textit{173}, 1--16, doi:
  10.1111/j.1365--246X.2007.03,700.x, 2008.

\bibitem[{\textit{Miranian}(2004)}]{Miranian2004}
Miranian, L., Slepian functions on the sphere, generalized {G}aussian
  quadrature rule, \textit{Inv.~Prob.}, \textit{20}, 877--892, 2004.

\bibitem[{\textit{Mitra and Maniar}(2006)}]{Mitra+2006}
Mitra, P.~P., and H.~Maniar, Concentration maximization and local basis
  expansions ({LBEX}) for linear inverse problems, \textit{IEEE Trans.~Biomed
  Eng.}, \textit{53}(9), 1775--1782, 2006.

\bibitem[{\textit{Moore and Cada}(2004)}]{Moore+2004}
Moore, I.~C., and M.~Cada, Prolate spheroidal wave functions, an introduction
  to the {S}lepian series and its properties,
  \textit{Appl.~Comput.~Harmon.~Anal.}, \textit{16}, 208--230, 2004.

\bibitem[{\textit{Narcowich and Ward}(1996)}]{Narcowich+96}
Narcowich, F.~J., and J.~D. Ward, Nonstationary wavelets on the m-sphere for
  scattered data, \textit{Appl.~Comput.~Harmon.~Anal.}, \textit{3}, 324--336,
  1996.

\bibitem[{\textit{Nashed and Walter}(1991)}]{Nashed+91}
Nashed, M.~Z., and G.~G. Walter, General sampling theorems for functions in
  {R}eproducing {K}ernel {H}ilbert {S}paces, \textit{Math.~Control Signals
  Syst.}, \textit{4}, 363--390, 1991.

\bibitem[{\textit{Nystr\"om}(1930)}]{Nystrom30}
Nystr\"om, E.~J., {\"U}ber die praktische {A}ufl\"osung von
  {I}ntegralgleichungen mit {A}nwendungen auf {R}andwert\-aufgaben,
  \textit{Acta Mathematica}, \textit{54}, 185--204, 1930.

\bibitem[{\textit{Olhede and Walden}(2002)}]{Olhede+2002}
Olhede, S., and A.~T. Walden, Generalized {M}orse wavelets, \textit{IEEE
  Trans.~Signal~Process.}, \textit{50}(11), 2661--2670, 2002.

\bibitem[{\textit{Olhede and Metikas}(2009)}]{Olhede+2009}
Olhede, S.~C., and G.~Metikas, The monogenic wavelet transform, \textit{IEEE
  Trans.~Signal~Process.}, \textit{57}(9), 3426--3441, doi:
  10.1109/TSP.2009.2023,397, 2009.

\bibitem[{\textit{Panet et~al.}(2006)\textit{Panet, Chambodut, Diament,
  Holschneider, and Jamet}}]{Panet+2006}
Panet, I., A.~Chambodut, M.~Diament, M.~Holschneider, and O.~Jamet, New
  insights on intraplate volcanism in {F}rench {P}olynesia from wavelet
  analysis of {GRACE}, {CHAMP}, and sea surface data,
  \textit{J.~Geophys.~Res.}, \textit{111}, B09,403, doi: 10.1029/2005JB004,141,
  2006.

\bibitem[{\textit{Papoulis}(1975)}]{Papoulis75}
Papoulis, A., A new algorithm in spectral analysis and band-limited
  extrapolation, \textit{IEEE-CS}, \textit{22}(9), 735--742, 1975.

\bibitem[{\textit{Parks and Shenoy}(1990)}]{Parks+90}
Parks, T.~W., and R.~G. Shenoy, Time-frequency concentrated basis functions, in
  \textit{Proc.~IEEE Int.~Conf.~Acoust.~Speech Signal Process.}, vol.~5, pp.
  2459--2462, IEEE, 1990.

\bibitem[{\textit{Parlett and Wu}(1984)}]{Parlett+84}
Parlett, B.~N., and W.-D. Wu, Eigenvector matrices of symmetric tridiagonals,
  \textit{Numer.~Math.}, \textit{44}, 103--110, 1984.

\bibitem[{\textit{Percival and Walden}(1993)}]{Percival+93}
Percival, D.~B., and A.~T. Walden, \textit{Spectral Analysis {f}or Physical
  Applications, Multitaper {a}nd Conventional Univariate Techniques}, Cambridge
  Univ.~Press, New York, 1993.

\bibitem[{\textit{Percival and Walden}(2006)}]{Percival+2006}
Percival, D.~B., and A.~T. Walden, \textit{Wavelet methods for time series
  analysis}, Cambridge Univ.~Press, 2006.

\bibitem[{\textit{Press et~al.}(1992)\textit{Press, Teukolsky, Vetterling, and
  Flannery}}]{Press+92}
Press, W.~H., S.~A. Teukolsky, W.~T. Vetterling, and B.~P. Flannery,
  \textit{Numerical Recipes {i}n {FORTRAN}: {T}he Art {o}f Scientific
  Computing}, 2nd ed., Cambridge Univ.~Press, New York, 1992.

\bibitem[{\textit{Ramesh and Lean}(1991)}]{Ramesh+91}
Ramesh, P.~S., and M.~H. Lean, Accurate integration of singular kernels in
  boundary integral formulations for {H}elmholtz equation,
  \textit{Int.~J.~Num.~Meth.~Eng.}, \textit{31}, 1055--1068, 1991.

\bibitem[{\textit{Riedel and Sidorenko}(1995)}]{Riedel+95}
Riedel, K.~S., and A.~Sidorenko, Minimum bias multiple taper spectral
  estimation, \textit{IEEE Trans.~Signal~Process.}, \textit{43}(1), 188--195,
  1995.

\bibitem[{\textit{Saito}(2007)}]{Saito2007}
Saito, N., Data analysis and representation on a general domain using
  eigenfunctions of {L}aplacian, \textit{Appl.~Comput.~Harmon.~Anal.},
  \textit{25}, 68--97, doi: 10.1016/j.acha.2007.09.005, 2007.

\bibitem[{\textit{Schmidt et~al.}(2006)\textit{Schmidt, Han, Kusche, Sanchez,
  and Shum}}]{Schmidt+2006}
Schmidt, M., S.-C. Han, J.~Kusche, L.~Sanchez, and C.~K. Shum, Regional
  high-resolution spatiotemporal gravity modeling from {GRACE} data using
  spherical wavelets, \textit{Geophys.~Res.~Lett.}, \textit{33}(8), L0840, doi:
  10.1029/2005GL025,509, 2006.

\bibitem[{\textit{Schmidt et~al.}(2007)\textit{Schmidt, Fengler,
  Mayer-G{\"u}rr, Eicker, Kusche, S{\'a}nchez, and Han}}]{Schmidt+2007}
Schmidt, M., M.~Fengler, T.~Mayer-G{\"u}rr, A.~Eicker, J.~Kusche,
  L.~S{\'a}nchez, and S.-C. Han, Regional gravity modeling in terms of
  spherical base functions, \textit{J.~Geodesy}, \textit{81}(1), 17--38, doi:
  10.1007/s00,190--006--0101--5, 2007.

\bibitem[{\textit{Schott and Th\'ebault}(2011)}]{Schott+2011}
Schott, J.-J., and E.~Th\'ebault, Modelling the earth’s magnetic field from
  global to regional scales, in \textit{Geomagnetic Observations {a}nd Models},
  \textit{IAGA Special Sopron Book Ser.}, vol.~5, edited by M.~Mandea and
  M.~Korte, Springer, 2011.

\bibitem[{\textit{Shepp and Zhang}(2000)}]{Shepp+2000}
Shepp, L., and C.-H. Zhang, Fast functional magnetic resonance imaging via
  prolate wavelets, \textit{Appl.~Comput.~Harmon.~Anal.}, \textit{9}(2),
  99--119, doi: 10.1006/acha.2000.0302, 2000.

\bibitem[{\textit{Shkolnisky}(2007)}]{Shkolnisky2007}
Shkolnisky, Y., Prolate spheroidal wave functions on a disc --- {I}ntegration
  and approximation of two-dimensional bandlimited functions,
  \textit{Appl.~Comput.~Harmon.~Anal.}, \textit{22}, 235--256, doi:
  10.1016/j.acha.2006.07.002, 2007.

\bibitem[{\textit{Shkolnisky et~al.}(2006)\textit{Shkolnisky, Tygert, and
  Rokhlin}}]{Shkolnisky+2006}
Shkolnisky, Y., M.~Tygert, and V.~Rokhlin, Approximation of bandlimited
  functions, \textit{Appl.~Comput.~Harmon.~Anal.}, \textit{21}, 413--420, doi:
  10.1016/j.acha.2006.05.001, 2006.

\bibitem[{\textit{Simons}(2010)}]{Simons2010}
Simons, F.~J., Slepian functions and their use in signal estimation and
  spectral analysis, in \textit{Handbook of Geomathematics}, edited by
  W.~Freeden, M.~Z. Nashed, and T.~Sonar, chap.~30, pp. 891--923, doi:
  10.1007/978--3--642--01,546--5\_30, Springer, Heidelberg, Germany, 2010.

\bibitem[{\textit{Simons and Dahlen}(2006)}]{Simons+2006b}
Simons, F.~J., and F.~A. Dahlen, Spherical {S}lepian functions and the polar
  gap in geodesy, \textit{Geophys.~J.~Int.}, \textit{166}, 1039--1061, doi:
  10.1111/j.1365--246X.2006.03,065.x, 2006.

\bibitem[{\textit{Simons and Dahlen}(2007)}]{Simons+2007}
Simons, F.~J., and F.~A. Dahlen, A spatiospectral localization approach to
  estimating potential fields on the surface of a sphere from noisy, incomplete
  data taken at satellite altitudes, in \textit{Wavelets {XII}}, vol. 6701,
  edited by D.~{Van de Ville}, V.~K. Goyal, and M.~Papadakis, pp. 670,117, doi:
  10.1117/12.732,406, SPIE, 2007.

\bibitem[{\textit{Simons et~al.}(2003)\textit{Simons, van~der Hilst, and
  Zuber}}]{Simons+2003a}
Simons, F.~J., R.~D. van~der Hilst, and M.~T. Zuber, Spatio-spectral
  localization of isostatic coherence anisotropy in {A}ustralia and its
  relation to seismic anisotropy: {I}mplications for lithospheric deformation,
  \textit{J.~Geophys.~Res.}, \textit{108}(B5), 2250, doi:
  10.1029/2001JB000,704, 2003.

\bibitem[{\textit{Simons et~al.}(2006)\textit{Simons, Dahlen, and
  Wieczorek}}]{Simons+2006a}
Simons, F.~J., F.~A. Dahlen, and M.~A. Wieczorek, Spatiospectral concentration
  on a sphere, \textit{SIAM Rev.}, \textit{48}(3), 504--536, doi:
  10.1137/S0036144504445,765, 2006.

\bibitem[{\textit{Simons et~al.}(2009)\textit{Simons, Hawthorne, and
  Beggan}}]{Simons+2009b}
Simons, F.~J., J.~C. Hawthorne, and C.~D. Beggan, Efficient analysis and
  representation of geophysical processes using localized spherical basis
  functions, in \textit{Wavelets~{XIII}}, vol. 7446, edited by V.~K. Goyal,
  M.~Papadakis, and D.~{Van de Ville}, pp. 74,460G, doi: 10.1117/12.825,730,
  SPIE, 2009.

\bibitem[{\textit{Simons et~al.}(1997)\textit{Simons, Solomon, and
  Hager}}]{Simons+97a}
Simons, M., S.~C. Solomon, and B.~H. Hager, Localization of gravity and
  topography: Constraints on the tectonics and mantle dynamics of {V}enus,
  \textit{Geophys.~J.~Int.}, \textit{131}, 24--44, 1997.

\bibitem[{\textit{Slepian}(1964)}]{Slepian64}
Slepian, D., Prolate spheroidal wave functions, {F}ourier analysis and
  uncertainty --- {IV}: {E}xtensions to many dimensions; generalized prolate
  spheroidal functions, \textit{Bell Syst.~Tech.~J.}, \textit{43}(6),
  3009--3057, 1964.

\bibitem[{\textit{Slepian}(1976)}]{Slepian76}
Slepian, D., On bandwidth, \textit{Proc.~IEEE}, \textit{64}(3), 292--300, 1976.

\bibitem[{\textit{Slepian}(1978)}]{Slepian78}
Slepian, D., Prolate spheroidal wave functions, {F}ourier analysis and
  uncertainty --- {V}: {T}he discrete case, \textit{Bell Syst.~Tech.~J.},
  \textit{57}, 1371--1429, 1978.

\bibitem[{\textit{Slepian}(1983)}]{Slepian83}
Slepian, D., Some comments on {F}ourier analysis, uncertainty and modeling,
  \textit{SIAM Rev.}, \textit{25}(3), 379--393, 1983.

\bibitem[{\textit{Slepian and Pollak}(1961)}]{Slepian+61}
Slepian, D., and H.~O. Pollak, Prolate spheroidal wave functions, {F}ourier
  analysis and uncertainty --- {I}, \textit{Bell Syst.~Tech.~J.},
  \textit{40}(1), 43--63, 1961.

\bibitem[{\textit{Slepian and Sonnenblick}(1965)}]{Slepian+65}
Slepian, D., and E.~Sonnenblick, Eigenvalues associated with prolate spheroidal
  wave functions of zero order, \textit{Bell Syst.~Tech.~J.}, \textit{44}(8),
  1745--1759, 1965.

\bibitem[{\textit{Tegmark}(1995)}]{Tegmark95}
Tegmark, M., A method for extracting maximum resolution power spectra from
  galaxy surveys, \textit{Astroph.~J.}, \textit{455}, 429--438, 1995.

\bibitem[{\textit{Tegmark}(1996)}]{Tegmark96a}
Tegmark, M., A method for extracting maximum resolution power spectra from
  microwave sky maps, \textit{Mon.~Not.~R.~Astron.~Soc}, \textit{280},
  299--308, 1996.

\bibitem[{\textit{Thomson}(1982)}]{Thomson82}
Thomson, D.~J., Spectrum estimation and harmonic analysis, \textit{Proc.~IEEE},
  \textit{70}(9), 1055--1096, 1982.

\bibitem[{\textit{Thomson}(1990)}]{Thomson90}
Thomson, D.~J., Quadratic-inverse spectrum estimates: applications to
  paleoclimatology, \textit{Phil.~Trans.~R.~Soc.~London, Ser.~A},
  \textit{332}(1627), 539--597, 1990.

\bibitem[{\textit{Tricomi}(1970)}]{Tricomi70}
Tricomi, F.~G., \textit{Integral Equations}, 5 ed., Interscience, New York,
  1970.

\bibitem[{\textit{van~de Ville and Unser}(2008)}]{Vandeville+2008}
van~de Ville, D., and M.~Unser, Complex wavelet bases, steerability, and the
  {M}arr-like pyramid, \textit{IEEE Trans.~Image~Proc.}, \textit{17}(11),
  2063--2080, doi: 10.1109/TIP.2008.2004,797, 2008.

\bibitem[{\textit{van~de Ville et~al.}(2002)\textit{van~de Ville, Philips, and
  Lemahieu}}]{Vandeville+2002}
van~de Ville, D., W.~Philips, and I.~Lemahieu, On the {N}-dimensional extension
  of the discrete prolate spheroidal window, \textit{IEEE
  Trans.~Signal~Process.}, \textit{9}(3), 89--91, 2002.

\bibitem[{\textit{Walden}(1990)}]{Walden90a}
Walden, A.~T., Improved low-frequency decay estimation using the multitaper
  spectral-analysis method, \textit{Geophys.~Prospect.}, \textit{38}, 61--86,
  1990.

\bibitem[{\textit{Walter and Soleski}(2005)}]{Walter+2005b}
Walter, G., and T.~Soleski, A new friendly method of computing prolate
  spheroidal wave functions and wavelets, \textit{Appl.~Comput.~Harmon.~Anal.},
  \textit{19}, 432--443, 2005.

\bibitem[{\textit{Walter and Shen}(2004)}]{Walter+2004}
Walter, G.~G., and X.~Shen, Wavelets based on prolate spheroidal wave
  functions, \textit{J.~Fourier Anal.~Appl.}, \textit{10}(1), 1--26,
  doi:10.1007/s00,041--004--8001--7, 2004.

\bibitem[{\textit{Walter and Shen}(2005)}]{Walter+2005a}
Walter, G.~G., and X.~Shen, Wavelet like behavior of {S}lepian functions and
  their use in density estimation, \textit{Comm.~Stat.~Theory Meth.},
  \textit{34}(3), 687--711, 2005.

\bibitem[{\textit{Walter and Soleski}(2008)}]{Walter+2008}
Walter, G.~G., and T.~Soleski, Error estimates for the {PSWF} method in {MRI},
  \textit{Contemp.~Math.}, \textit{451}, 262, 2008.

\bibitem[{\textit{Wei et~al.}(2010)\textit{Wei, Kennedy, and
  Lamahewa}}]{Wei+2010}
Wei, L., R.~A. Kennedy, and T.~A. Lamahewa, Signal concentration on unit
  sphere: {A}n azimuthally moment weighting approach, in \textit{Proc.~IEEE
  Int.~Conf.~Acoust.~Speech Signal Process.}, pp. 1--4, IEEE, 2010.

\bibitem[{\textit{Wieczorek and Simons}(2005)}]{Wieczorek+2005}
Wieczorek, M.~A., and F.~J. Simons, Localized spectral analysis on the sphere,
  \textit{Geophys.~J.~Int.}, \textit{162}(3), 655--675, doi:
  10.1111/j.1365--246X.2005.02,687.x, 2005.

\bibitem[{\textit{Wieczorek and Simons}(2007)}]{Wieczorek+2007}
Wieczorek, M.~A., and F.~J. Simons, Minimum-variance spectral analysis on the
  sphere, \textit{J.~Fourier Anal.~Appl.}, \textit{13}(6), 665--692,
  10.1007/s00,041--006--6904--1, 2007.

\bibitem[{\textit{Wingham}(1992)}]{Wingham92}
Wingham, D.~J., The reconstruction of a band-limited function and its {F}ourier
  transform from a finite number of samples at arbitrary locations by
  {S}ingular {V}alue {D}ecomposition, \textit{IEEE Trans.~Signal~Process.},
  \textit{40}(3), 559--570, doi: 10.1109/78.120,799, 1992.

\bibitem[{\textit{Xiao et~al.}(2001)\textit{Xiao, Rokhlin, and
  Yarvin}}]{Xiao+2001}
Xiao, H., V.~Rokhlin, and N.~Yarvin, Prolate spheroidal wavefunctions,
  quadrature and interpolation, \textit{Inverse Problems}, \textit{17},
  805--838, 10.1088/0266--5611/17/4/315, 2001.

\bibitem[{\textit{Yang et~al.}(2002)\textit{Yang, Lindquist, Shepp, Zhang,
  Wang, and Smith}}]{Yang+2002}
Yang, Q.~X., M.~A. Lindquist, L.~Shepp, C.-H. Zhang, J.~Wang, and M.~B. Smith,
  Two dimensional prolate spheroidal wave functions for {MRI},
  \textit{J.~Magnet.~Reson.}, \textit{158}, 43--51, 2002.

\bibitem[{\textit{Yao}(1967)}]{Yao67}
Yao, K., Application of reproducing kernel {H}ilbert spaces --- {B}andlimited
  signal models, \textit{Inform.~Control}, \textit{11}(4), 429--444, 1967.

\bibitem[{\textit{Zhang}(1994)}]{Zhang94}
Zhang, X., Wavenumber spectrum of very short wind waves: {A}n application of
  two-dimensional {S}lepian windows to spectral estimation,
  \textit{J.~Atm.~Ocean.~Tech.}, \textit{11}, 489--505, 1994.

\bibitem[{\textit{Zhou et~al.}(1984)\textit{Zhou, Rushforth, and
  Frost}}]{Zhou+84}
Zhou, Y., C.~K. Rushforth, and R.~L. Frost, Singular value decomposition,
  singular vectors, and the discrete prolate spheroidal sequences,
  \textit{Proc.~IEEE Int.~Conf.~Acoust.~Speech Signal Process.}, \textit{9}(1),
  92--95, 1984.

\end{thebibliography}
\bibliographystyle{agufull}

\end{document}